\documentclass[a4paper, draft]{article}
\usepackage[all]{xy}
\usepackage{amssymb,amsmath,amsthm}
\numberwithin{equation}{section}

\newcommand{\Q}{\mathbb{Q}}

\newcommand{\R}{\mathbb{R}}
\newcommand{\C}{\mathbb{C}}
\newcommand{\N}{\mathbb{N}}
\newcommand{\Z}{\mathbb{Z}}

\newtheorem{thm}{Theorem}
\newtheorem{lem}{Lemma}

\theoremstyle{definition}
\newtheorem*{ack}{Acknowledgements}

\renewcommand{\mod}[1]{\hspace{-2.9mm}\pmod{#1}}
\newcommand{\x}{{\bf x}}

\newcommand{\e}{\emph}
\newcommand{\rom}{\mathrm}
\newcommand{\bfP}{\mathbb{P}}
\newcommand{\A}{\mathbb{A}}
\newcommand{\ov}{\overline}
\newcommand{\ma}{\mathbf}
\newcommand{\ben}{\begin{enumerate}}
\newcommand{\een}{\end{enumerate}}
\newcommand{\eit}{\begin{itemize}}
\newcommand{\beq}{\begin{equation}}
\newcommand{\eeq}{\end{equation}}

\newcommand{\ve}{\varepsilon}

\newcommand{\mcal}{\mathcal}
\newcommand{\lab}{\label}
\newcommand{\al}{\alpha}
\newcommand{\om}{\omega}

\newcommand{\D}{\Delta}
\newcommand{\be}{\beta}

\newcommand{\hcf}{\rom{gcd}}
\newcommand{\colt}[2]{\genfrac{}{}{0pt}{1}{#1}{#2}}
\newcommand{\tolt}[3]{\colt{#1}{\colt{#2}{#3}}}
\newcommand{\eqm}[2]{\equiv #1 \pmod{#2}}

\newcommand{\Dfive}{{\mathbf D}_5}
\renewcommand{\leq}{\leqslant}
\renewcommand{\geq}{\geqslant}
\renewcommand{\d}{\mathrm{d}}
\newcommand{\vr}{\varrho}
\newcommand{\vt}{\vartheta}
\renewcommand{\d}{\mathrm{d}}

\begin{document}

\title{On Manin's conjecture for singular del Pezzo surfaces of degree
  four, I}

\author{R. de la Bret\`eche$^1$ and T.D. Browning$^2$\\
\small{$^1$\emph{Universit\'e Paris-Sud,
B\^atiment 425, 91405 Orsay Cedex}}\\
\small{$^{2}$\emph{School of Mathematics, Bristol University, Bristol
    BS8 1TW}}\\ 
\small{$^1$regis.de-la-breteche@math.u-psud.fr,
  $^2$t.d.browning@bristol.ac.uk}}
\date{}

\maketitle

\begin{abstract}
This paper contains a proof of the Manin conjecture for the singular
del Pezzo surface
$$
X: \quad x_0x_1-x_2^2 = x_0x_4-x_1x_2+x_3^2 = 0,
$$
of degree four.
In fact, if $U \subset X$ is the open subset formed by deleting the
unique line 
from $X$, and $H$ is the usual height function on $\bfP^4(\Q)$, then
the height zeta function 
$
\sum_{x \in U(\Q)}{H(x)^{-s}}
$
is analytically continued to the half-plane $\Re e (s)>9/10$. 
\end{abstract}

\section{Introduction}

Let $Q_1,Q_2 \in \Z[x_0,\ldots,x_4]$ be a pair of quadratic forms whose common zero
locus defines a geometrically integral surface $X\subset \bfP^4$.
Then $X$ is a del Pezzo surface of degree four.  We assume
henceforth that the set
$X(\Q)=X\cap \bfP^4(\Q)$ of rational points on $X$ is non-empty,
so that in particular $X(\Q)$ is dense in $X$ under
the Zariski topology.
Given a point $x=[x_0,\ldots,x_4] \in \bfP^4(\Q),$ with
$x_0,\ldots,x_4 \in \Z$ such that $\hcf(x_0,\ldots,x_4)=1$, we let
$
H(x)=\max_{0\leq i \leq 4}|x_i|.
$
Then $H: \bfP^4(\Q) \rightarrow \R_{\geq 0}$ is
the height attached to the anticanonical divisor $-K_X$ on $X$,
metrized by the choice of norm $\max_{0\leq i \leq 4}|x_i|$.
A finer notion of density is  provided by analysing the
asymptotic behaviour of the quantity
$$
N_{U,H}(B)=\#\{x \in U(\Q): H(x) \leq B\},
$$
as $B\rightarrow \infty$,
for appropriate open subsets $U\subseteq X$.  Since every quartic del Pezzo
surface $X$ contains a line, it is natural to estimate $N_{U,H}(B)$ for the
open subset $U$ obtained by deleting the lines from $X$.  
The motivation behind  this paper is to consider the asymptotic
behaviour of $N_{U,H}(B)$ for singular del Pezzo surfaces of
degree four.

A classification of quartic del Pezzo surfaces $X \subset \bfP^4$ can
be found in the work of Hodge and Pedoe \cite[Book IV, \S
XIII.11]{h-p}, showing in particular that there are only finitely
many isomorphism classes to consider.
Let $\tilde X$ denote the minimal desingularisation of $X$, and let
$\rom{Pic}{\,\tilde X}$ be the Picard group of  $\tilde X$.
Then Manin has stated a very general conjecture \cite{f-m-t} that
predicts the asymptotic behaviour of counting functions
associated to suitable Fano varieties.
In our setting it leads us to expect
the existence of a positive constant $c_{X,H}$ such that
\begin{equation}\lab{manin}
N_{U,H}(B) = c_{X,H} B (\log B)^{\rho-1}\big(1+o(1)\big),
\end{equation}
as $B \rightarrow \infty,$
where $\rho$ denotes the rank of $\rom{Pic}{\,\tilde X}$.
The constant $c_{X,H}$ has received a conjectural interpretation
at the hands of Peyre \cite{p}, which in turn
has been generalised to cover certain other cases by Batyrev and
Tschinkel~\cite{b-t} and 
Salberger~\cite{salb}. A brief discussion of $c_{X,H}$ will take place in \S
\ref{conform}.

There has been rather little progress towards the Manin conjecture for del Pezzo
surfaces of degree four.  The main successes in this direction are to
be found in work of Batyrev and Tschinkel  \cite{b-t'}, covering the
case of toric varieties, and in the work of Chambert-Loir and
Tschinkel \cite{ct}, covering the case of equivariant
compactifications of vector groups.
It is our intention to investigate the distribution of rational points
in the special case that $X$ is defined by the pair of equations
\begin{align*}
x_0x_1-x_2^2&=0,\\
x_0x_4-x_1x_2+x_3^2 &= 0.
\end{align*}
Then $X \subset \bfP^4$ is a del Pezzo surface of degree four, with  
a unique singular point $[0,0,0,0,1]$ which is of type $\Dfive$.
Furthermore, $X$ contains precisely one line
$
x_0=x_2=x_3=0.
$
It turns out that $X$ is an equivariant compactification of
$\mathbb{G}_a^2$, so that the work of Chambert-Loir and
Tschinkel \cite[Theorem 0.1]{ct} ensures that the asymptotic formula
\eqref{manin} holds when $U\subset X$ is taken to 
be the open subset formed by deleting the unique line from $X$.
Nonetheless there are several reasons why this problem is still
worthy of attention.  Firstly, in making an exhaustive study of $X$
it is hoped that a template will 
be set down for the treatment of other singular del Pezzo
surfaces.   In fact no explicit use is made of the fact that $X$ is an equivariant compactification of
$\mathbb{G}_a^2$, and the techniques that we develop in this paper 
have already been applied to other surfaces \cite{2,3}.
Secondly, in addition to improving upon Chambert-Loir and Tschinkel's asymptotic formula for $N_{U,H}(B)$,
the results that we obtain lend themselves more readily as a
bench-test for future refinements of the Manin conjecture, such as
that recently proposed by Swinnerton-Dyer \cite{swd} for example.

Let $X\subset \bfP^4$ be the $\Dfive$ del Pezzo
surface defined above, and let $U\subset X$ be the corresponding open subset.
Our first result concerns the height zeta function
\begin{equation}\lab{h-zeta}
Z_{U,H}(s)=\sum_{x \in U(\Q)}\frac{1}{H(x)^s},
\end{equation}
that is defined when $\Re e(s)$ is sufficiently large.   
The analytic properties of $Z_{U,H}(s)$ are intimately related to
the asymptotic behaviour of the counting function $N_{U,H}(B)$.  
For $\Re e (s)>0$ we define the functions
\begin{align}
E_1(s+1)&= \zeta(6s+1)\zeta(5s+1)\zeta(4s+1)^2\zeta(3s+1)
\zeta(2s+1),\lab{e1}\\
E_2(s+1)&=
\frac{ \zeta(14s+3)\zeta(13s+3)^3}{\zeta(10s+2)\zeta(9s+2)
\zeta(8s+2)^3\zeta(7s+2)^3\zeta(19s+4)}.\lab{e2}
\end{align}
It is easily seen that $E_1(s)$ has a meromorphic 
continuation to the entire complex plane with a single pole at $s=1$.
Similarly it is clear that $E_2(s)$ is holomorphic and bounded 
on the half-plane $\Re e (s)\geq 9/10+\ve$, for any $\ve>0$. 
We are now ready to state our main result.

\begin{thm}\lab{main'}
Let $\ve>0$.  
Then there exists a constant $\beta \in \R$, 
and functions $G_1(s), G_2(s)$
that are holomorphic   on the half-plane $\Re
e(s)\geq 5/6+\ve$, such that for $\Re e(s)>1$ we have
$$
Z_{U,H}(s) =E_1(s)E_2(s)G_1(s) +
\frac{12/\pi^2+2\beta}{s-1}+G_2(s). 
$$
In particular $(s-1)^6Z_{U,H}(s)$ has a holomorphic 
continuation to the half-plane $\Re e(s)> 9/10.$ 
The function $G_1(s)$ is bounded for $\Re e(s)\geq 5/6+\ve$ and satisfies $G_1(1)\neq 0$, 
and the function $G_2(s)$ satisfies  
$$ 
G_2(s)\ll_\ve (1+|\Im m (s)|)^{6(1-\Re e (s))+\varepsilon} 
$$
on this domain.
\end{thm}

An explicit expression for $\beta$ can be found in
\eqref{beta'}, whereas the formulae \eqref{G1}--\eqref{saw}  
can be used to deduce an explicit expression for $G_1$.
There are several features of Theorem \ref{main'} that are worthy of remark.  The first
step in the  proof of Theorem \ref{main'} is the observation
\begin{equation}\lab{integral}
Z_{U,H}(s)=s\int_1^\infty t^{-s-1}N_{U,H}(t)\d t.
\end{equation}
Thus we find ourselves in the situation of establishing a preliminary
estimate for $N_{U,H}(B)$ in order to deduce the analytic properties
of $Z_{U,H}(s)$ presented in Theorem \ref{main'}, before then using
this information to deduce an improved estimate for $N_{U,H}(B)$. This
will be given in Theorem \ref{main} below. 
With this order of events in mind we highlight that the term 
$\frac{12}{\pi^2}(s-1)^{-1}$ appearing in Theorem~\ref{main'} 
corresponds to an isolated conic contained in $X$. 
Moreover, whereas the first term $E_1(s)E_2(s)G_1(s)$ in the expression for
$Z_{U,H}(s)$ corresponds to the main term in our preliminary estimate
for $N_{U,H}(B)$, and 
arises through the 
approximation of certain arithmetic quantities by real-valued
continuous functions, the term involving $\beta$ has a more arithmetic
interpretation.  Indeed, it will be seen to arise purely out of
the error terms produced by approximating these arithmetic quantities
by continuous functions. 
Finally we make the observation that under the assumption of the
Riemann hypothesis $E_2(s)$ is holomorphic for $\Re
e(s)> 17/20$, so that $Z_{U,H}(s)$ has an analytic continuation to
this domain.

%R%pour m\'emoire
%La majoration \eqref{G2} peut para\^{\i}tre d\'ecevante puisqu'elle est 
%la simple cons\'e\-quence de $G_2(s)\ll 1+|\Im m (s)|$ dans tout demi-plan
%de la forme 
%$\{ s\in
%\C  :\Re e(s)\geq 5/6+\ve\}$ via le th\'eor\`eme de Phragmen-Lindel\"of.
%Cependant elle est toujours plus forte que toutes les majorations
%ponctuelles de $Z_{U,H}(s)$ que l'on peut esp\'erer du fait de l'\'etat
%des connaissances actuelles sur la fonction $\zeta$ (voir la
%majoration \eqref{majE1}).

We now come to how Theorem \ref{main'} can be used to deduce an asymptotic formula for
$N_{U,H}(B)$.
We shall verify in \S \ref{conform} that the following result is in
accordance with the Manin conjecture.

\begin{thm}\lab{main}
Let $\delta \in (0,1/12)$.
Then there exists a polynomial $P$ of
degree $5$ such that 
$$
N_{U,H}(B)=B P(\log B) +O(B^{1-\delta}),
$$
for any $B \geq 1$.
Moreover the leading coefficient of $P$ is equal to
$$
\frac{\tau_\infty}{28800}
\prod_{p}\Big(1-\frac{1}{p}\Big)^{6 }
\Big(1+\frac{6}{p}+
\frac{1}{p^2}\Big),
$$
where
\begin{equation}\lab{t-inf}
\tau_\infty = \int_{0}^1
\int_{-1}^{1/v}\Big(\min\{\sqrt{u^3+1},1/v^3\}-
\sqrt{\max\{u^3-1,0\}}\Big)\d u\d v.
\end{equation}
\end{thm}

The deduction of Theorem \ref{main} from Theorem \ref{main'} will take
place in \S \ref{deduc}, and amounts to a routine application of Perron's formula.  
Although we choose not to give the details here, 
it is in fact possible to take  $\delta \in
(0,1/11)$ in the statement of Theorem \ref{main}
by using more sophisticated estimates for moments of the Riemann zeta
function in the critical strip.
By expanding the height zeta function $Z_{U,H}(s)$ as
a power series in $(s-1)^{-1}$, one may obtain
explicit expressions for the lower order coefficients of the polynomial $P$
in Theorem \ref{main}.
It would be interesting to obtain refinements of Manin's conjecture
that admit conjectural interpretations of the lower order coefficients.

The principal tool in the proof of Theorem \ref{main'} is a passage to the universal
torsor above the minimal desingularisation $\tilde X$ of $X$.
Although originally introduced to aid in the study of the Hasse principle and Weak approximation, universal torsors
have recently become endemic in the context of counting rational 
points of bounded height.  In \S \ref{prelim} we shall
establish a bijection between $U(\Q)$ and integer points on the universal torsor, 
which in this setting has the natural affine embedding 
$$
v_2y_0^2y_4 -  v_0y_1^3y_2^2 + v_3y_3^2=0.
$$
Once we have translated the problem to the universal torsor,
Theorem \ref{main'} will be established using a range of techniques drawn from 
classical analytic number theory. In particular Theorem \ref{main'} seems to
be the first instance of a height zeta function that has been calculated via 
a passage to the universal torsor.

\begin{ack}
Parts of this paper were worked upon while the authors were attending
the conference \e{Diophantine geometry} at G\"ottingen University in
June 2004, during which time the authors benefited from useful   
conversations with several participants.  
While working on this paper, the first author was at the \'Ecole
Normale Sup\'erieure,  and the second author was at Oxford University supported by 
EPSRC grant number GR/R93155/01.
The hospitality and financial support of these institutions is
gratefully acknowledged.
Finally, the authors would like to thank the anonymous
referee for his careful reading of the manuscript.
\end{ack}

\section{Conformity with the Manin conjecture}\lab{conform}

In this section we shall show that Theorem \ref{main} agrees with the
Manin conjecture.  For this we need to calculate the value of
$c_{U,H}$ and $\rho$ in \eqref{manin}.  We therefore
review some of the geometry of the surface $X\subset
\bfP^4$, as defined by the pair of quadratic forms
\begin{equation}\lab{forms}
Q_1(\x)=x_0x_1-x_2^2, \quad Q_2(\x)=x_0x_4-x_1x_2+x_3^2,
\end{equation}
where $\x=(x_0,\ldots,x_4)$.

Let $\tilde X$ denote the minimal desingularisation
of $X$,  and let $\pi: \tilde{X}\rightarrow X$ denote the corresponding
blow-up map.
We let $E_6$ denote the strict transform of the unique line
contained in $X$, and let $E_1,\ldots,E_5$ denote the exceptional
curves of~$\pi$.  Then the divisors $E_1,\ldots,E_6$
satisfy the Dynkin diagram
$$
\xymatrix{
       & & E_3 \ar@{-}[d] \\
E_5 \ar@{-}[r] & E_4 \ar@{-}[r] & E_1 \ar@{-}[r] & E_2
\ar@{-}[r] & E_6 }
$$
and generate the Picard group of $\tilde{X}$. In particular we have
$\rho=6$ in \eqref{manin}, which agrees with Theorem \ref{main}.

It remains to discuss the conjectured value of the constant $c_{X,H}$ in
\eqref{manin}.  For this we shall follow the presentation of
Batyrev and Tschinkel \cite[\S 3.4]{b-t}.
Let $\Lambda_\rom{eff}(\tilde{X})\subset \rom{Pic}{\,\tilde X} \otimes_\Z
\R$ be the cone of effective divisors on $\tilde X$, and let
$$
\Lambda_\rom{eff}^\vee(\tilde{X}) =\{\ma{s} \in \rom{Pic}^\vee{\tilde
X} \otimes_\Z
\R: \langle\ma{s}, \ma{t} \rangle \geq 0 \,\,\forall \,\ma{t} \in
\Lambda_\rom{eff}(\tilde{X})\}
$$
be the corresponding dual cone, where $\rom{Pic}^\vee{\tilde X}$
denotes the dual lattice to $\rom{Pic}{\,\tilde X}.$
Then if $\d\ma{t}$ denotes the Lebesgue measure on
$\rom{Pic}^\vee{ \tilde X} \otimes_\Z \R$, we define
$$
\alpha(\tilde X) =
\int_{\Lambda_\rom{eff}^\vee(\tilde{X})} e^{-\langle-K_{\tilde
X},\ma{t} \rangle}\d\ma{t},
$$
where $-K_{\tilde X}$ is the anticanonical divisor of $\tilde
X$.
Thus $\alpha(\tilde X)$ measures the volume of the polytope obtained by
intersecting $\Lambda_\rom{eff}^\vee(\tilde{X})$ with a certain
affine hyperplane.

Next we discuss the Tamagawa measure on the closure
$\ov{\tilde{X}(\Q)}$ of $\tilde{X}(\Q)$ in $\tilde{X}(\A_\Q)$, where
$\A_\Q$ denotes the adele ring. Write
$L_p(s,\rom{Pic}{\,\tilde{X}})$ for the local factors of
$L(s,\rom{Pic}{\,\tilde{X}})$.  Furthermore, let
$\omega_\infty$ denote the archimedean density of points on $X$, and let
$\omega_p$ denote the usual $p$-adic density of points on $X$, for any prime $p$.  Then we
may define the Tamagawa measure
\begin{equation}\lab{formulatauH}
\tau_H(\tilde{X})= \lim_{s\rightarrow 1}\big((s-1)^\rho
L(s,\rom{Pic}{\,\tilde{X}})\big)
\omega_\infty\prod_p
\frac{\omega_p}{L_p(1,\rom{Pic}{\,\tilde{X}})},
\end{equation}
where $\rho$ denotes the rank of $\rom{Pic}{\,\tilde X}$ as above.
With these definitions in mind, the
conjectured value of the constant in \eqref{manin} is equal to
\begin{equation}\lab{constant}
c_{X,H}=\alpha(\tilde X) \beta(\tilde X) \tau_H(\tilde X),
\end{equation}
where $\beta(\tilde X)=\#
H^1(\mathrm{Gal}(\overline{\Q}/\Q),\rom{Pic}{\,\tilde X} \otimes_\Q
\overline{\Q})=1$, since $\tilde X$ is split over the ground field $\Q$.

We begin
by calculating the value of $\alpha(\tilde X)$, for which we shall
follow the approach of Peyre and Tschinkel \cite[\S 5]{p-t}.  We
need to determine the cone $\Lambda_\rom{eff}(\tilde{X})$ and the
anticanonical divisor $-K_{\tilde X}$.  To determine these we may
use the Dynkin diagram to write down the intersection matrix
\begin{center}
\begin{tabular}{r|rrrrrr}
     & $E_1$ & $E_2$ & $E_3$ & $E_4$ & $E_5$ & $E_6$\\
\hline
$E_1$ & $-2$ & $1$ & $1$ & $1$ & $0$ & $0$\\
$E_2$ & $1$ & $-2$ & $0$ & $0$ & $0$ & $1$\\
$E_3$ & $1$ & $0$ & $-2$ & $0$ & $0$ & $0$\\
$E_4$ & $1$ & $0$ & $0$ & $-2$ & $1$ & $0$\\
$E_5$ & $0$ & $0$ & $0$ & $1$ & $-2$ & $0$\\
$E_6$ & $0$ & $1$ & $0$ & $0$ & $0$ & $-1$
\end{tabular}.
\end{center}
But then it follows that
$\Lambda_\rom{eff}(\tilde{X})$ is generated by the divisors $E_1,\ldots,E_6$.
Furthermore, on employing the adjunction formula $-K_{\tilde X}.D=2+D^2$ for $D
\in \rom{Pic}{\,\tilde X}$, we easily deduce that
$$
-K_{\tilde X}=6E_1+5E_2+3E_3+4E_4+2E_5+4E_6.
$$
It therefore follows that
\begin{align*}
\alpha(\tilde X)&= \rom{Vol}\big\{(t_1,t_2,t_3,t_4,t_5,t_6)\in
    \R_{\geq 0}^6:
6t_1+5t_2+3t_3+4t_4+2t_5+4t_6=1 \big\}\\
&= \frac{1}{4}
\rom{Vol}\big\{(t_1,t_2,t_3,t_4,t_5)\in \R_{\geq 0}^5:
6t_1+5t_2+3t_3+4t_4+2t_5\leq 1
\big\},
\end{align*}
whence in fact
\begin{equation}\lab{alpha}
\alpha(\tilde X)=1/345600.
\end{equation}
It remains to calculate the value of
$\tau_H(\tilde X)$ in \eqref{constant}.

\begin{lem}\lab{localdensities}
We have $\tau_H(\tilde{X})= 12\tau_{\infty}\tau$, where $\tau_\infty$
is given by \eqref{t-inf} and 
\begin{equation}\lab{deftauf}
\tau=\prod_{ p
}\Big(1-\frac{1}{p}\Big)^{6 }\Big(1+\frac{6}{p}+
\frac{1}{p^2}\Big).
\end{equation}
\end{lem}

\begin{proof}
Recall the definition \eqref{formulatauH} of $\tau_H(\tilde{X})$.
Our starting point is the observation that
$
L(s,\rom{Pic}{\,\tilde{X}})=\zeta(s)^6.
$
Hence it easily follows that
$$
\lim_{s\rightarrow 1}\big((s-1)^\rho
L(s,\rom{Pic}{\,\tilde{X}})\big)=
    \lim_{s\rightarrow 1}\big((s-1)^6 \zeta(s)^6\big)
=1.
$$
Furthermore we plainly  have
\begin{equation}\lab{Lp1}
L_p(1,\rom{Pic}{\,\tilde{X}})^{-1}=
\Big(1-\frac{1}{p}\Big)^6,
\end{equation}
for any prime $p$.

We proceed by employing the method of Peyre \cite{p} to calculate the value of the
archimedean density $\omega_\infty$.
It will be convenient to parameterise
the points via the choice of variables $x_0,x_1,x_4$,
for which we first observe that the Leray form 
$\omega_L(\tilde X)$ is given by $(4x_2x_3)^{-1}\d x_0\d x_1\d x_4,$ since
$$
\det \left(
\begin{matrix}
\frac{\partial Q_1}{ \partial x_2}&\frac{\partial Q_2}{
\partial x_2} \cr\cr  \frac{\partial Q_1}{ \partial x_3}
&\frac{\partial Q_2}{ \partial x_3}
\end{matrix}
\right)=-4x_2x_3.
$$
Now in any real solution to the pair of equations $Q_1(\x)=Q_2(\x)=0$,
the components $x_0$ and $x_1$ must necessarily share the same
sign.  Taking into account the fact that $\x$ and $-\x$ represent the
same point in $\bfP^4$, we therefore see that
$$
\omega_{\infty}  = 2
\int_{ \{\x\in\R^5: ~
Q_1(\x)=Q_2(\x)=0,~
0\leq x_0,x_1,x_3, |x_4| \leq 1\}}\omega_L( \tilde X).
$$
We write $\omega_{\infty,-}$ to denote the contribution to
$\omega_\infty$ from the case $x_2=-\sqrt{x_0x_1}$, and
$\omega_{\infty,+}$ for the contribution from the case $x_2=\sqrt{x_0x_1}$.
But then it follows that
$$
\omega_{\infty,+} =\frac{1}{2} \int\int\int
\frac{\d x_0\d x_1\d x_4}{\sqrt{x_0^2x_1(x_0^{-1/2}x_1^{3/2}-x_4)}},
$$
where the triple integral is over all $x_0,x_1,x_4\in \R$ such that
$$
0\leq x_0,x_1,|x_4| \leq 1,  \quad x_1^{3/2}\geq
x_0^{1/2}x_4, \quad 1+x_0x_4\geq x_1^{3/2}x_0^{1/2}.
$$
The change of variables $u=x_1^{1/2}/ x_0^{1/6}$, and then
$v=x_0^{1/6}$, therefore yields
$$
\omega_{\infty,+} =6\int\int\int
\frac{\d u\d v\d x_4}{\sqrt{u^3-x_4}},
$$
where the triple integral is now over all $u,v,x_4\in \R$ such that
$$
0\leq v\leq 1, \quad 0\leq u\leq 1/v, \quad \max\{-1,u^3-1/v^6\}\leq x_4
\leq \min\{1,u^3\}.
$$
On performing the integration over $x_4$, a straightforward
calculation leads to the equality
$$
\omega_{\infty,+}=12\int_{0}^1
\int_{0}^{1/v}\Big(\min\{\sqrt{u^3+1},1/v^3\}-
\sqrt{\max\{u^3-1,0\}}\Big)\d u\d v.
$$
The calculation of $\omega_{\infty,-}$ is similar.  On
following the steps outlined above, one is easily led to the equality
%% $$
%% \omega_{\infty,-} =-6\int\int\int
%% \frac{\d u\d v\d x_4}{\sqrt{-u^3-x_4}},
%% $$
%% where the triple integral is over all $u,v,x_4\in \R$ such that
%% $$
%% 0\leq v\leq 1, \quad 0\leq u\leq 1/v, \quad \max\{-1,-u^3-1/v^6\}\leq x_4
%% \leq -u^3.
%% $$
%% A simple integration over $x_4$, followed by an obvious change of
%% variables,  therefore yields
$$
\omega_{\infty,-}=12\int_{0}^1
\int_{-1}^{0}\min\{\sqrt{u^3+1},1/v^3\}\d u\d v.
$$
Once taken together, these equations combine to show that $\om_\infty=12 \tau_\infty$, where
$\tau_\infty$ is given by \eqref{t-inf}.

It remains to calculate the value of
$\omega_p=\lim_{r\to \infty}p^{-3r}N(p^r)$, where we have written
$N(p^r)=\#\{\x ~(\bmod{p^r}): Q_1(\x)\equiv
Q_2(\x)\equiv 0 ~(\bmod{p^r})\}.$
To begin with we write
$$
x_0=p^{k_0}x_0',\qquad x_1=p^{k_1}x_1'
$$
with $p\nmid x_0'x_1'$.
Now we have $p^r \mid x_2^2$ if and only if $k_0+k_1\geq r$,
and there are at most $p^{r/2}$ square roots of zero modulo $p^r$.
When $ k_0+k_1< r$, it follows that $k_0+k_1$ must be even and we may write
$$
x_2=p^{(k_0+k_1)/2}x_2',
$$
with $p\nmid x_2'$ and
$$
x_0'x_1'-{x_2'}^2\equiv 0 \mod{p^{r-k_0-k_1}}.
$$
The number of possible choices for $x_0',x_1',x_2'$ is therefore
$$
h_p(r,k_0,k_1)=\left\{
\begin{array}{ll} \phi(p^{r-k_0})\phi(p^{r-(k_0+k_1)/2})p^{k_0} &
\mbox{if $k_0+k_1< r$},\\
O(p^{5r/2-k_0- k_1 })& \mbox{if $k_0+k_1\geq r$}.
\end{array}
\right.
$$

It remains to determine the number of solutions $x_3,x_4$ modulo $p^r$
such that
\begin{equation}\lab{equationx3x4}
p^{k_0}{x_0'}x_4-p^{k_0/2+3k_1/2}x_1'x_2'+x_3^2\equiv 0 \mod{p^r}.
\end{equation}
In order to do so we distinguish between three basic cases:  either
$k_0+k_1<r$ and $k_0\leq 3k_1$, or
$k_0+k_1<r$ and $k_0>3k_1$, or
else $k_0+k_1\geq r$. For the first
two of these cases we must take care only to sum over values of
$k_0,k_1$ such that $k_0+k_1$ is even.
We shall denote by $N_i(p^r)$ the contribution to $N(p^r)$ from the
$i$th case, for $1\leq i \leq 3$, so that
\begin{equation}\label{N_i}
N(p^r)=N_1(p^r)+N_2(p^r)+N_3(p^r).
\end{equation}

We begin by calculating the value of $N_1(p^r)$.  For this we write
$x_3=p^{k_3}x_3'$, with $k_3=\min \{r/2, \lceil k_0/2 \rceil\}=\lceil
k_0/2 \rceil$.   The number of
possibilities for $x_3'$ is $p^{r-\lceil k_0/2 \rceil}$, each one
leading to a congruence of the form
$$
{x_0'}x_4- p^{3k_1/2-k_0/2}x_1'x_2'  +p^{2\lceil k_0/2
    \rceil-k_0}{x_3'}^2 \equiv 0 \mod{p^{r-k_0}}.
$$
Modulo $p^{r-k_0}$, there is one choice for $x_4$,
and so there are $p^{r+k_0-\lceil k_0/2 \rceil}=
p^{r+\lfloor k_0/2 \rfloor}$ possibilities for
$x_3$ and $x_4$.
On summing these contributions over all the relevant
values of $k_0,k_1$, we therefore obtain
$$
N_1(p^r)= \sum_{\colt{k_0+k_1<r, ~k_0,k_1\geq 0}
{k_0\leq 3k_1, ~2\mid (k_0+k_1)}}
p^{r+\lfloor k_0/2 \rfloor} h_p(r,k_0,k_1)
= p^{3r-2}(p^2+4p+1)\big(1+o(1)\big),
$$
as $r\rightarrow \infty$.

Next we calculate $N_2(p^r)$, for which we shall not use the above
  calculation for $h_p(r,k_0,k_1)$.  On writing
$x_3=p^{k_3}x_3'$, with
$$k_3=\min \{r/2, \lceil k_0/4+3k_1/4 \rceil\}=
\lceil k_0/4+3k_1/4 \rceil,
$$
we observe that $k_0/2+3k_1/2$ must be even since $p \nmid x_1'x_2'$.
Thus $k_3=k_0/4+3k_1/4$ and
$p\nmid x_3'$.  In this way \eqref{equationx3x4} becomes
$$
p^{k_0/2-3k_1/2}x_0'x_4 -x_1'x_2'+{x_3'}^2\equiv 0 \mod{p^{r-k_0/2-3k_1/2}},
$$
which thereby implies that ${x_3'}^2\equiv x_1'x_2' ~~\mod{p^{k_0/2-3k_1/2}}$.
At this point we recall the auxiliary congruence
${x_2'}^2\equiv x_0'x_1' ~~\mod{p^{r-k_0-k_1}}$ that is
satisfied by $x_0',x_1'x_2'$.
We proceed by fixing values of $x_2'$ and $x_3'$, for which there are
precisely $(1-1/p)^2p^{r-k_0/2-k_1/2}p^{r-k_3}$ choices.
But then $x_1'$ is fixed modulo ${p^{k_0/2-3k_1/2}}$, and so there are
$p^{r-k_1-(k_0/2-3k_1/2)}$ possibilities for $x_1'$.
Finally we deduce from the remaining two congruences that
there are  $p^{k_1}$ ways of fixing $x_0'$, and $p^{k_0}$ ways
of fixing $x_4$.
Summing over the relevant values
of $k_0$ and $k_1$ we therefore obtain
$$
N_2(p^r)
=\sum_{\colt{k_0+k_1<r, ~k_0,k_1\geq 0}
{k_0>3k_1, ~4\mid (k_0-k_1)}}
(1-1/p)^2p^{3r-(k_0-k_1)/4} = 2p^{3r-1}\big(1+o(1)\big),
$$
as $r\rightarrow \infty$.

Finally we calculate the value of $N_3(p^r)$.  In this case we
write $x_2=p^{r/2}x_2'$.  But then a similar 
calculation ultimately shows that $N_3(p^r)=o(p^{3r})$ 
as $r\rightarrow \infty$.
On combining our estimates for $N_1(p^r),N_2(p^r), N_3(p^r)$ into
\eqref{N_i}, we
therefore deduce that
$$
N(p^r)=p^{3r}\left(1+\frac{6}{p}+\frac{1}{p^2}\right)\big(1+o(1)\big)
$$
as $r \rightarrow \infty$, whence
$$
\omega_p=\lim_{r\to\infty}p^{-3r}N(p^r)=1+\frac{6}{p}+\frac{1}{p^2}
$$
for any prime $p$. We combine this with \eqref{Lp1}, in the manner indicated by
\eqref{formulatauH}, in order to deduce \eqref{deftauf}.
\end{proof}

We end this section by combining \eqref{alpha} and Lemma
\ref{localdensities} in \eqref{constant}, in order to deduce that the
conjectured value of the constant in \eqref{manin} is equal to
$$
c_{X,H}=\frac{1}{28800}\tau_{\infty}\tau,
$$
where $\tau_\infty$ is given by \eqref{t-inf} and
$\tau$ is given by \eqref{deftauf}.
This agrees with the value of the leading coefficient obtained in Theorem \ref{main}.

\section{Congruences}\lab{lin-cong}

In this section we shall collect together some of the basic
facts concerning congruences that will be needed in the proof of 
Theorem \ref{main'}.  We begin by discussing the case of quadratic congruences.
For any integers $a,q$ such that $q>0$, we
define the arithmetic function
$\eta(a;q)$ to be the number of positive integers $n\leq q$
such that
$
n^2\eqm{a}{q}.
$
When $q$ is odd it follows that 
$$
\eta(a;q)=\sum_{d \mid q}|\mu(d)| \big(\frac{a}{d}\big),
$$
where $(\frac{a}{d})$ is the usual Jacobi
symbol. On noting that $\eta(a;2^{\nu})\leq 4$ 
for any $\nu\in \N$, it easily follows that
\begin{equation}\lab{triv}
\eta(a;q) \leq  2^{\omega(q)+1} 
\end{equation}
for any $q \in \N$.  Here $\omega(q)$
denotes the number of distinct prime factors of $q$.

Turning to the case of linear congruences, let
$\kappa\in[0,1]$ and let 
$\vt$ be any arithmetic function such that
\begin{equation}\lab{condtheta}
\sum_{d=1}^\infty \frac{|(\vartheta*\mu)(d) |}{d^\kappa}<\infty,
\end{equation}
where $(f*g)(d)=\sum_{e\mid d}f(e)g(d/e)$ is the usual
Dirichlet convolution of any two arithmetic functions $f,g$.
Then for any coprime integers $a,q$ such that $q>0$, and any
$t\geq 1$, we deduce that
\begin{align*}
\sum_{\colt{n\leq t}{n\equiv
a\mod{q}}}\vartheta(n)
&=\sum_{\colt{d=1}{\hcf(d,q)=1}}^\infty(\vartheta*\mu)(d) \sum_{\colt{m\leq
t/d}{md\equiv a\mod{q}}}1\nonumber\\
&=\frac{t}{q}\sum_{\colt{d=1}{\hcf(d,q)=1}}^{\infty}
\frac{(\vartheta*\mu)(d)}{d}
+O\Big(t^\kappa\sum_{d=1}^\infty \frac{|(\vartheta*\mu)(d) |}{d^\kappa}\Big),
\end{align*}
on using the equality $\vt=(\vt*\mu)*1$
and the trivial estimate $\lfloor x \rfloor=x+O(x^{\kappa})$ for any $x>0$. We
summarise this estimate in the following result.

\begin{lem}\lab{tek}
Let $\kappa\in[0,1]$, let $\vt$ be any arithmetic function such that
\eqref{condtheta} holds, and let $a,q \in \Z$ be such that $q>0$ and
$\hcf(a,q)=1$. Then we have
$$
\sum_{\colt{n\leq t}{n\equiv
a\mod{q}}}\vartheta(n)
=\frac{t}{q}\sum_{\colt{d=1}{\hcf(d,q)=1}}^{\infty}
\frac{(\vartheta*\mu)(d)}{d}
+O\Big(t^\kappa\sum_{d=1}^\infty \frac{|(\vartheta*\mu)(d) |}{d^\kappa}\Big).
$$
\end{lem}

Define the real-valued function $\psi(t)=\{t\}-1/2$, where
$\{t\}$ denotes the fractional part of $t \in \R$.  Then
$\psi$ is periodic with period $1$.
When $\vartheta(n)=1$ for all $n \in \N$ we are able to 
refine Lemma \ref{tek} considerably.

\begin{lem}\lab{cong}
Let $a, q \in \Z$ be such that $q>0$, and let $t_1, t_2\in \R$ such
that $t_2 \geq \max\{0,t_1\}$.
  Then
$$
\#\{t_1<n \leq t_2: n\equiv a\mod{q} \}=\frac{t_2-t_1}{q}+ r(t_1,t_2;a,q),
$$
where
$$
r(t_1,t_2;a,q)=\psi\Big( \frac{t_1-a}{q}\Big) - \psi\Big( \frac{t_2-a}{q}\Big).
$$
\end{lem}

\begin{proof}
Write $a=b+qc$ for some integer $0\leq b<q$. 
Then it is clear that
$$
\#\{t_1< n \leq t_2: n\equiv 
b\mod{q}\}=\Big\lfloor \frac{t_2-b}{q}\Big\rfloor
-\Big\lfloor \frac{t_1-b}{q}\Big\rfloor,
$$
whence 
$$
\#\{t_1<n \leq t_2: n\equiv a\mod{q}\}-\frac{t_2-t_1}{ q}= r(t_1,t_2;b,q).
$$
We complete the proof of Lemma \ref{cong} by  noting that
$r(t_1,t_2;b,q)=r(t_1,t_2;a,q)$, since $\psi$ has  period $1$.
\end{proof}

We shall also need to know something about the average order of the
function $\psi$.
We proceed by demonstrating the following result.

\begin{lem}\lab{av-order}
Let $\ve>0$, $t \geq 0$ and let $X \geq 1$.  Then for any $b,q \in \Z$ such that
$q>0$ and $\hcf(b,q)=1$, we have
$$
\sum_{ 0\leq x<X}
\psi\Big(\frac{t-bx^2}{q}\Big)\ll_\varepsilon 
(qX)^{\varepsilon}\Big(\frac{X}{q^{1/2}}+q^{1/2}\Big).
$$
\end{lem}

\begin{proof}
Throughout this proof we shall write $e(t)=e^{2\pi it}$ and 
$e_q(t)=e^{2\pi it/q}$.
In order to establish Lemma \ref{av-order} 
we shall expand the function
$f(k)=\psi( {(t-k)/q} )$ as a Fourier series. Thus we have
$$
f(k)=\sum_{0\leq \ell<q}a(\ell)e_q(k\ell),
$$
for any $k\in\Z$, where the coefficients $a(\ell)$ are given by
$$
a(\ell)=\frac{1}{q}\sum_{0\leq j<q}f(j)e_q(-j\ell).
$$
Let $\|\alpha\|$ denote the distance from $\al \in \R$ to the nearest integer.
We proceed by proving the estimates
\begin{equation}\lab{majcoeff}
a(\ell)\ll \left\{
\begin{array}{ll}
q^{-1}, & \ell=0,\\
q^{-1}\parallel \ell/q\parallel^{-1}, & \ell \neq 0.
\end{array}
\right.
\end{equation}
This is straightforward.  To verify the estimate for $a(0)$ we
simply note that
\begin{align*}
a(0)
&=\frac{1}{q}\sum_{0\leq j<q} \Big(\Big\{\frac{t-j}{q}\Big\}-\frac{1}{2}\Big)\\
&=\frac{1}{q}\sum_{0\leq j\leq t} \Big(\frac{t-j-q/2}{q}\Big) 
+\frac{1}{q}\sum_{t<
j<q}\Big(\frac{t-j+q/2}{q}\Big)
\ll\frac{1}{ q}.
\end{align*}
Similarly, when $\ell\neq 0$ we have
\begin{align*}
a(\ell)&=\frac{1}{q}\sum_{0\leq j\leq t}
\Big(\frac{t-j-q/2}{q}\Big)e_q(-j\ell)+
\frac{1}{q}\sum_{t<j<q}\Big(\frac{t-j+q/2}{q}\Big)e_q(-j\ell)\\
&=\frac{1}{q}\sum_{0\leq j<q} \frac{-j }{q}e_q(-j\ell)-\frac{1}{2q}\sum_{0\leq
j\leq t} e_q(-j\ell)+\frac{1}{2q}\sum_{t<j< q} e_q(-j\ell)
\ll\frac{1}{q\parallel \ell/q\parallel},
\end{align*}
as required.

In view of the above we therefore obtain
\begin{align*}
\sum_{ 0\leq x<X} \psi\Big(\frac{t-bx^2}{q}\Big)
&=\sum_{0\leq \ell<q}a(\ell)\sum_{ 0\leq x<X}e_q( \ell bx^2)\\
&= a(0)\lceil X \rceil +\sum_{m\mid q}\sum_{\colt{1\leq \ell'<q/m}{\hcf 
(\ell',q/m)=1}} 
\hspace{-0.3cm}
a(\ell'm)\sum_{ 0\leq
x<X}e_{q/m}(\ell'x^2).
\end{align*}
But here the inner sum can plainly be estimated using Weyl's
inequality, and so has size
$$
\ll_\varepsilon  X^\ve\Big(\frac{m^{1/2}X}{q^{1/2}} + \frac{q^{1/2}}{m^{1/2}}\Big).
$$
On employing
\eqref{majcoeff}, we therefore deduce that
\begin{align*}
\sum_{ 0\leq x<X}
\psi\Big(\frac{t-bx^2}{q}\Big)&\ll_\ve
\frac{X}{q}+  \sum_{m\mid q}m^{1/2}X^\ve\sum_{1\leq
  \ell'<q/m}\frac{ {q^{-1/2}}X+ q^{1/2}}{q
   \parallel \ell'm/q\parallel}\\
&\ll_\ve (qX)^{2\varepsilon}\Big(\frac{X}{q^{1/2}}+q^{1/2}\Big),
\end{align*}
which completes the proof of Lemma \ref{av-order}.
\end{proof}

Let $\ve>0$ and let $t \geq 0$.  Then for any $b,q \in \Z$ such that
$q>0$ and $\hcf(b,q)=1$, we may deduce from Lemma \ref{av-order} that
\begin{equation}\lab{MAJ}
\sum_{ 0\leq x<q}
\psi\Big(\frac{t-bx^2}{q}\Big)\ll_\varepsilon q^{1/2+\varepsilon}.
\end{equation}
But then it follows from an application of 
M\"obius inversion that
\begin{align*}
\sum_{\colt{0\leq x<q}{\hcf(x,q)=1}}
\psi\Big(\frac{t-bx^2}{q}\Big)
&=\sum_{n \mid q } \mu (n )\sum_{ 0\leq x'<q/n }
\psi\Big(\frac{t/n  -b n  {x'}^2}{q /n}\Big)\\
&=\sum_{\colt{n \mid q }{m=\hcf(n,q/n)}} 
\hspace{-0.3cm}
\mu (n )m
\hspace{-0.3cm}
\sum_{ 0\leq x'<q/(mn)}
\hspace{-0.3cm}
\psi\Big(\frac{t/(mn)  -bn{x'}^2/m}{q/(mn)}\Big),
\end{align*}
whence \eqref{MAJ} yields
$$
\sum_{\colt{0\leq x<q}{\hcf(x,q)=1}}
\psi\Big(\frac{t-bx^2}{q}\Big)
\ll_\varepsilon q^{1/2+\varepsilon}
\sum_{n \mid q} \Big(\frac{\hcf(n,q/n)}{n}\Big)^{1/2} \ll_\ve q^{1/2+
2\varepsilon}.
$$
This therefore establishes the following result, 
on re-defining the choice of $\ve$.

\begin{lem}\lab{av-order'}
Let $\ve>0$ and let $t \geq 0$.  Then for any $b,q \in \Z$ such that
$q>0$ and $\hcf(b,q)=1$, we have
$$
\sum_{\colt{0\leq x<q}{\hcf(x,q)=1}}
\psi\Big(\frac{t-bx^2}{q}\Big) \ll_\varepsilon q^{1/2+ \varepsilon}.
$$
\end{lem}

\section{Preliminary manoeuvres}\lab{prelim}

We begin this section by introducing some notation.
For any $n \geq 2$ we let $Z^{n+1}$ denote the set of
primitive vectors in $\Z^{n+1}$,
where $\ma{v}=(v_0,\ldots,v_n) \in \Z^{n+1}$ is
said to be primitive if
$
\hcf(v_0,\ldots,v_n)=1.
$
Moreover we shall let $\Z_*^{n+1}$ (resp. $Z_*^{n+1}$) denote the
set of vectors $\ma{v} \in \Z^{n+1}$ (resp. $\ma{v} \in Z^{n+1}$) such
that $v_0\cdots v_n \neq 0$.  Finally we underline the fact that
throughout our work $\N$ is always taken to denote the set of 
positive integers.

The proof of Theorem \ref{main'} rests
upon establishing a preliminary asymptotic formula for the counting
function $N_{U,H}(B)$.
Recall the definition \eqref{forms} of the quadratic forms
$Q_1,Q_2$.
Our first task in this section is to relate $N_{U,H}(B)$ to the quantity
$$
N(Q_1,Q_2;B)=\#\{\x\in Z_*^5:  0<x_0 ,x_1,x_3\leq B,|x_4|\leq B,
~Q_1(\x)=Q_2(\x)=0\}.
$$
In fact we shall establish the following result rather easily.

\begin{lem}\lab{red-1}
Let $B\geq 1$.  Then we have
$$
N_{U,H}(B) =2N(Q_1,Q_2;B) +\frac{12}{\pi^2}B+O(B^{2/3}).
$$
\end{lem}

\begin{proof}
It is clear that any solution to the pair of equations
$Q_1(\x)=Q_2(\x)=0$ which satisfies $x_0=0$, must in fact correspond to a
point lying on the line $x_0=x_2=x_3=0$ contained in $X$.  On noting
that $\x$ and $-\x$ represent the same point in
projective space, we therefore deduce that
$$
N_{U,H}(B)=\frac{1}{2}\#\{\x\in Z^5: \|\x\|\leq B,
~Q_1(\x)=Q_2(\x)=0, ~x_0 \neq 0\},
$$
where $\|\x\|=\max_{0\leq i \leq 4}|x_i|$.
We proceed to consider the contribution from the vectors $\x \in
Z^5$ for which $\|\x\|\leq B$ and
\begin{equation}\lab{test}
Q_1(\x)=Q_2(\x)=0, \quad x_1x_2x_3x_4=0.
\end{equation}
Note first that $x_1=0$ if and only if $x_2=0$ in \eqref{test}, since
$x_0 \neq 0$.
Thus if we consider the contribution from those vectors for which
$x_1x_2=0$, it follows that we must count integers $|x_0|,|x_3|,|x_4| 
\leq B$ for which
$\hcf(x_0,x_3,x_4)=1$ and $x_0x_4+x_3^2=0$.
Now either $x_3=0$, in which case $\x=(1,0,0,0,0)$ since $\x$ is
primitive and $x_0 \neq 0$, or else the primitivity of $\x$ implies that
$\x=(a^2,0,0,\pm ab,- b^2)$ 
for coprime non-zero integers $a, b$.  Hence the overall contribution
from this case is clearly  $12B/\pi^2+O(B^{1/2}).$

Suppose now that  $x_3=0$  and $x_1x_2\neq 0$ in \eqref{test}.
Then we must count the number of mutually coprime non-zero integers 
$x_0,x_1,x_2,x_4$,
with modulus at most $B$, such that $x_0x_1=x_2^2$ and  $x_0x_4=x_1x_2$.
Since we are only interested in an upper bound it clearly suffices
to count non-zero
integers $x_0,x_1, x_4$, with modulus at most $B$,  such that 
$\hcf(x_0,x_1,x_4)=1$ and $x_0x_4^2=x_1^3$.  
But then it follows that $(x_0,x_1,x_4)=\pm (a^3, ab^2, b^3)$  
for coprime non-zero integers $a,b$, 
whence the overall contribution is $O(B^{2/3})$.
Finally the case $x_4=0$ and $x_1x_2\neq 0$ in \eqref{test} is
handled in much the
same way, now via the parameterisation $\x=\pm (a^4,b^4, a^2b^2, ab^3,0)$.
This therefore
establishes that
$$
N_{U,H}(B)=\frac{1}{2}\#\{\x\in Z_*^5: \|\x\|\leq B,
~Q_1(\x)=Q_2(\x)=0\}+\frac{12}{\pi^2}B+O (B^{2/3}).
$$

We complete the proof of Lemma \ref{red-1} by choosing $x_0>0$ and
$x_3>0$. This then forces the inequality $x_1>0$, whence $\|\x\|=
\max\{ x_0,x_1,x_3,|x_4|\}.$
\end{proof}

We now turn to the task of establishing a bijection between the points
counted by $N(Q_1,Q_2;B)$ and integral points on the universal torsor
above the minimal desingularisation of $X$.
Let $\x \in Z_*^5$ be any vector counted by $N(Q_1,Q_2;B)$.
In particular it follows that $x_0,x_1,x_3$ are positive.
We begin by considering
solutions to the equation
$Q_1(\x)=0.$
But it is easy to see that there is a bijection
between the set of integers $x_0,x_1,x_2$ such that
$x_0,x_1>0$ and $x_0x_1=x_2^2$, and the set of $x_0,x_1,x_2$ such that
$$
x_0=z_0^2z_2, \quad x_1=z_1^2z_2, \quad x_2=z_0z_1z_2,
$$
for non-zero integers $z_0, z_1,z_2$ such that $z_0,z_2>0$ and
\begin{equation}\lab{capri-0}
\hcf(z_0,z_1)=1.
\end{equation}
We now substitute these values into the equation
$Q_2(\x)=0,$  in order to obtain
\begin{equation}\lab{mondeo-1}
x_4z_0^2z_2-z_0z_1^3z_2^2+x_3^2=0.
\end{equation}
It is clear that $z_0z_2$ divides $x_3^2$. Hence we write
$$
z_0=v_0v_3{y_0''}^2, \quad z_2=v_2v_3{y_2''}^2,
$$
for $v_0,v_2,v_3,y_0'',y_2'' \in \N $ such that the products $v_0v_3$
and $v_2v_3$ are square-free, with $\hcf(v_0,v_2)=1$.  In particular
the product $v_0v_2v_3$ is clearly square-free.
We easily deduce that $v_0v_2v_3y_0''y_2''$ must divide $x_3$, whence there
exists ${y_3''} \in \N$ such  that
$x_3=v_0v_2v_3y_0''y_2''{y_3''}.$
Combining the various coprimality conditions arising from
\eqref{capri-0} and the definitions of $v_0,v_2$ and $v_3$, we
therefore obtain
\begin{equation}\lab{capri-1}
|\mu(v_0v_2v_3)|=1, \quad \hcf(v_0v_3y_0'',z_1)=1,
\end{equation}
where $\mu(n)$ denotes the M\"obius function for any non-zero integer
$n$.
On making the appropriate substitutions into \eqref{mondeo-1}, we deduce that
\begin{equation}\lab{mondeo-2}
v_0v_3x_4{y_0''}^2 - v_2v_3{y_2''}^2z_1^3 +
v_0v_2{y_3''}^2=0.
\end{equation}

At this point it is convenient to deduce a further coprimality
condition which follows from the assumption made at the outset that
$\hcf(x_0,\ldots,x_4)=1$.
Recalling the various changes of variables
that we have made so far, it is easily checked that
$
\hcf(x_0,x_1,x_2,x_3)=v_2v_3y_2''\hcf(y_2'',v_0y_0''y_3'').
$
Hence we find that
\begin{equation}\lab{capri-2}
\hcf(v_2v_3y_2'',x_4)=1.
\end{equation}
Now it follows from \eqref{mondeo-2} that $v_0$ divides
$v_2v_3{y_2''}^2z_1^3$.  But then since $v_0$ is square-free, we may
conclude
from \eqref{capri-1} that $v_0 \mid{y_2''}$.  Similarly we deduce from
\eqref{capri-1} and  \eqref{capri-2} that
$v_2 \mid y_0''$ and $v_3 \mid {y_3''}$.  Thus there exist
$y_0', y_2',y_3'\in \N$ and $ y_1 ,y_4 \in \Z_*$
such that
$$
y_0''=v_2y_0', \quad z_1=y_1, \quad{y_2''}=v_0y_2', \quad {y_3''}=v_3y_3',
\quad x_4=y_4.
$$
Substituting the above into \eqref{mondeo-2} we therefore obtain the
equation
\begin{equation}\lab{mondeo-3}
v_2{y_0'}^2y_4 -  v_0y_1^3{y_2'}^2 + v_3{y_3'}^2=0.
\end{equation}
Moreover we may combine \eqref{capri-1} and \eqref{capri-2} to get
\begin{equation}\lab{capri-3}
|\mu(v_0v_2v_3)|=1, \quad
\hcf(v_0v_2v_3{y_0'},y_1)=\hcf(v_0v_2v_3y_2',y_4)=1.
\end{equation}
Finally we write $v_1=\hcf({y_0'},{y_2'},y_3')$.  Thus there exist
$y_0,y_2,y_3 \in \N$ such that
$$
{y_0'}=v_1y_0, \quad {y_2'}=v_1y_2, \quad y_3'=v_1y_3,
$$
and we obtain the final equation
\begin{equation}\lab{ut}
v_2y_0^2y_4 -  v_0y_1^3y_2^2 + v_3y_3^2=0.
\end{equation}

It remains to collect together the coprimality conditions that have
arisen from this last change of variables.  First however we take a
moment to deduce three further coprimality conditions
\begin{equation}\lab{got}
\hcf(y_0,y_2)=1, \quad \hcf(y_0,y_3)=1, \quad\hcf(y_2,y_3)=1.
\end{equation}
To do so we simply use the obvious fact that $\hcf(y_0,y_2,y_3)=1$.
Suppose that $p$ is any prime divisor of
$y_2$ and $y_3$.   Then we clearly have $p^2 \mid v_2y_0^2y_4$ in
\eqref{ut}.  This is impossible by \eqref{capri-3} and the fact that
$\hcf(y_0,y_2,y_3)=1$. From this we may establish the second relation
in \eqref{got}.  Indeed, if $p\mid y_0,y_3$ then clearly $p^2 \mid
v_0y_1^3y_2^2$, which is impossible by \eqref{capri-3} and the fact
that $\hcf(y_2,y_3)=1$.  One checks the first relation in \eqref{got}
in a similar fashion.
Combining \eqref{capri-3} with \eqref{got} we therefore deduce that
\begin{equation}\lab{capri-4'}
\hcf(y_3,y_0y_2)=\hcf(y_4,v_0v_1v_2v_3y_2)=1
\end{equation}
and
\begin{equation}\lab{capri-5'}
|\mu(v_0v_2v_3)|=1, \quad \hcf(y_1,v_0v_1v_2v_3y_0)=
\hcf(y_0,y_2)=1.
\end{equation}

In fact it will be necessary to reformulate these coprimality
conditions somewhat.  We claim that once taken together with
\eqref{ut}, the relations \eqref{capri-4'} and \eqref{capri-5'} are
equivalent to
\begin{equation}\lab{capri-4}
\hcf(y_3,v_0y_0y_2)=\hcf(y_4,v_1v_2)=1
\end{equation}
and
\begin{align}
&\hcf(y_1,v_0v_1v_2v_3y_0)=1,
\label{capri-5}\\
|\mu(v_0v_2v_3)|=1,&\quad \hcf(v_2v_3y_0,y_2)=\hcf(v_0v_3,y_0 )=1.
\label{temple}
\end{align}
We first show how \eqref{ut}, \eqref{capri-4'} and \eqref{capri-5'}
imply \eqref{ut}, \eqref{capri-4}, \eqref{capri-5} and 
\eqref{temple}.  Suppose that $p$ is any prime divisor of $v_0$ and $y_3$.  Then \eqref{ut}
implies that $p\mid v_2y_0^2y_4$ which is easily seen to be impossible
via \eqref{capri-4'} and \eqref{capri-5'}.
Thus $\hcf(y_3,v_0)=1$.  Now suppose that $p$ is a prime divisor of
$v_3$ and $y_2$.  Then $p\mid v_2y_0^2y_4$ which is also impossible, and so
$\hcf(v_3,y_2)=1$.
The supplementary conditions $\hcf(v_2,y_2)=\hcf(v_0v_3,y_0 )=1$
easily follow from the relations $\hcf(v_0
y_2,v_3y_3)=\hcf(v_3,y_1) =1$, in addition to \eqref{ut}.
The converse is established along similar lines.

At this point we may summarise our argument as follows.
Let $\mcal{T} \subset \Z_*^{9}$ denote the set of
$(\ma{v,y})=(v_0,v_1,v_2,v_3,y_0,\ldots,y_4) \in \N^4\times \Z_*^5$
such that $y_0,y_2,y_3>0$,
\eqref{ut},  and \eqref{capri-4}--\eqref{temple} hold.
Then for any $\x
\in Z_*^5$ counted by $N(Q_1,Q_2;B)$, we have shown that there exists
$(\ma{v,y}) \in
\mcal{T}$ such that
\begin{align*}
x_0 &= v_0^4v_1^6v_2^5v_3^3y_0^4y_2^2,\\
x_1 &= v_0^2v_1^2v_2v_3y_1^2y_2^2, \\
x_2 &= v_0^3v_1^4v_2^3v_3^2y_0^2y_1 y_2^2,\\
x_3 &= v_0^2v_1^3v_2^2v_3^2y_0y_2y_3,\\
x_4 &= y_4.
\end{align*}
Conversely, given any
$(\ma{v,y}) \in \mcal{T}$ the point $\x$ given above
will be a solution of the equations $Q_1(\x)=Q_2(\x)=0$, with
$\x \in Z_*^5$.  To see the
primitivity of $\x$ we first recall that
once taken together with \eqref{ut}, the coprimality relations
\eqref{capri-4}--\eqref{temple} are equivalent to \eqref{capri-4'} 
and \eqref{capri-5'}.
But then it follows that
$\hcf(x_0,x_1,x_2,x_3)$ divides $ v_0^2v_1^2v_2v_3y_2^2$.  Finally
an application of
\eqref{capri-4'}  and \eqref{capri-5'} yields
\begin{align*}
\hcf(x_0,\ldots,x_4)&=\hcf(\hcf(x_0,x_1,x_2,x_3),x_4)\leq
\hcf(v_0^2v_1^2v_2v_3y_2^2,y_4)=1,
\end{align*}
as claimed.
Let us define the function $\Psi: \R^9 \rightarrow \R_{\geq 0}$, given by
$$
\Psi(\ma{v,y})=\max\big\{
    v_0^4v_1^6v_2^5v_3^3y_0^4y_2^2,
    v_0^2v_1^2v_2v_3y_1^2y_2^2,
    v_0^2v_1^3v_2^2v_3^2y_0y_2y_3, |y_4|
\big\}.
$$
We have therefore established the following result.

\begin{lem}\lab{base}
Let $B \geq 1$.  Then we have
$$
N(Q_1,Q_2;B)=\#\big\{(\ma{v,y}) \in \mcal{T}:
     \Psi(\ma{v,y}) \leq B\big\}.
$$
\end{lem}

It will become clear in subsequent sections that the equation
\eqref{ut} is a crucial ingredient in our proof of Theorem \ref{main'}.
In fact \eqref{ut} is an affine embedding of the universal torsor
above the minimal desingularisation of $X$.  Thus Derenthal, in work to
appear, has established the isomorphism
$$
\rom{Cox}(\tilde{X})= \mathrm{Spec}(\Q[\ma{v},\ma{y}]/(v_2y_0^2y_4 -  v_0y_1^3y_2^2 + v_3y_3^2)),
$$
where $\rom{Cox}(\tilde{X})$ is the Cox ring of $\tilde X$.

\section{The final count}

In this section we estimate $N(Q_1,Q_2;B)$, which we shall then combine with
Lemma \ref{red-1} to provide an initial estimate for $N_{U,H}(B)$.
Before proceeding with this task, it will be helpful to first outline our strategy.
In view of \eqref{ut} it is clear that for any
$(\ma{v,y})\in{\mcal T}$, the inequality $|y_4|\leq B$ is equivalent to
\begin{equation}\lab{encadrement}
-Bv_2y_0^2 \leq v_3y_3^2 - v_0y_1^3y_2^2 \leq Bv_2y_0^2.
\end{equation}
We henceforth write $\Phi(\ma{v,y})$ to denote the condition obtained
by replacing the term $|y_4|$ by $|(v_3y_3^2 -
v_0y_1^3y_2^2)/(v_2y_0^2)|$ in the
definition of $\Psi(\ma{v,y})$.

The basic idea behind our method is simply to view the equation
\eqref{ut} as a  congruence
$$
v_3y_3^2 \eqm{v_0y_1^3y_2^2}{v_2y_0^2}.
$$
Since we will have $\hcf(v_0y_1^3y_2^2,v_2y_0^2)=1$ when $(\ma{v,y})
\in \mcal{T}$, by \eqref{ut}, \eqref{capri-4} and \eqref{capri-5},
there exists a unique positive integer $\vr\leq
v_2y_0^2$ such that
$$
\hcf(\vr, v_2y_0^2)=1, \quad
v_3\vr^2 \eqm{v_0y_1}{v_2y_0^2},
$$
and
$$
y_3 \eqm{\vr y_1y_2}{v_2y_0^2}.
$$
The fact that $y_3$ and $y_4$ satisfy the coprimality conditions
\eqref{capri-4} complicates matters slightly, and makes it necessary to
first carry out a M\"obius inversion.

Next we analyse the inequality $\Phi(\ma{v,y})\leq B$.
In  doing so it will be convenient to define the quantities
\begin{equation}\lab{range1}
V_1=\left(\frac{B}{v_0^4v_2^5v_3^3y_0^4y_2^2}\right)^{1/6}
\end{equation}
and
\begin{equation}\lab{range2}
Y_1=\left(\frac{Bv_2y_0^2}{v_0y_2^2}\right)^{1/3}, \quad
Y_2=\Big(\frac{B}{v_0^4v_1^6v_2^5v_3^3y_0^4}\Big)^{1/2}, \quad
Y_3=\left(\frac{Bv_2y_0^2}{v_3}\right)^{1/2}.
\end{equation}
Moreover, we shall need to define the real-valued functions
$$
f_-(u,v)= \sqrt{\max\{u^3-1,0\}}, \quad
f_+(u,v)=\min \big\{\sqrt{ u^3+1}, 1/v^3\big\},
$$
and
\begin{equation}\lab{f}
f(u,v)=f_+(u,v)-f_-(u,v).
\end{equation}
In view of the inequality
$v_0^2v_1^3v_2^2v_3^2y_0y_2y_3\leq B$ that is implied by
$\Phi(\ma{v,y})\leq B$, we plainly have $y_3\leq V_1^3Y_3/v_1^3.$
A little thought therefore reveals that once combined with the
inequalities in \eqref{encadrement}, we have
\begin{equation}\lab{rolly}
Y_3 f_-(y_1/Y_1,v_1/V_1) \leq  y_3 \leq
Y_3 f_+(y_1/Y_1,v_1/V_1).
\end{equation}
Using the inequality $v_0^2v_1^2v_2v_3y_1^2y_2^2\leq B$, and 
deducing from \eqref{encadrement} that $y_1> -
Y_1$, we also see that
\begin{equation}\lab{holly}
-Y_1< y_1 \leq \frac{V_1Y_1}{v_1}.
\end{equation}
Next it follows from the inequality $\Phi(\ma{v,y})\leq B$ that
\begin{equation}\lab{bolly}
v_0^4v_1^6v_2^5v_3^3y_0^4y_2^2\leq B,
\end{equation}
whence
\begin{equation}\lab{golly}
1 \leq y_2 \leq Y_2
\end{equation}
and $1\leq v_1\leq V_1$.
In particular we must have $V_1\geq 1$, and so we may deduce the 
further inequality
\begin{equation}\lab{deV1}
V_1Y_1\leq V_1^3Y_1=\frac{B^{5/6}}{v_0^{7/3}
v_2^{13/6}v_3^{3/2}y_0^{4/3}y_2^{5/3}}.
\end{equation}
This will turn out to be useful at the end of \S \ref{y3y4}.

After having taken care of the contribution $S$, say, from the variables $y_3$
and $y_4$ in \S \ref{y3y4}, we will proceed in \S \ref{y1}
by summing $S$ over non-zero integers $y_1$ such that \eqref{holly}
holds and positive integers $y_2$ such that \eqref{golly} holds,
subject to certain conditions.  We shall denote this contribution by $S'$.
Finally, in \S \ref{final}, we shall obtain an estimate for
$N_{U,H}(B)$ by summing $S'$ over the remaining values of
$v_0,v_1,v_2,v_3,y_0$, subject to certain constraints, and then
applying Lemma \ref{red-1}.
During the course of the ensuing argument, in which we  establish
estimates for $S,S'$ and finally $N_{U,H}(B)$,
it will be convenient to handle the overall contribution from the 
error term in each estimate as we go.

\subsection{Summation over the variables $y_3$ and $y_4$}\lab{y3y4}

We begin by summing over the variables
$y_3,y_4$.
Let $(\ma{v},y_0,y_1,y_2) \in \N^4\times\Z_*^3$  satisfy
\eqref{capri-5}, \eqref{temple} and be
constrained to lie in the region defined by the inequalities
\eqref{holly}, \eqref{bolly}  and $y_0,y_2>0$.
As indicated above, we shall denote the double summation over $y_3$
and $y_4$ by $S$. In order to take care of the coprimality condition
$\hcf(y_4, v_1v_2)=1$ in \eqref{capri-4}, we apply a M\"obius
inversion to get
$$
S=\sum_{ k_4 \mid v_1v_2} \mu(k_4)S_{k_4},
$$
where the definition of $S_{k_4}$ is as for $S$ but with the extra
condition $k_4\mid y_4$ and without the
coprimality condition $\hcf(y_4, v_1v_2)=1$.
Thus it follows that  $S_{k_4}$ is equal to the number of non-zero
integers $y_3$ contained in the
region \eqref{rolly}, such that   $\hcf(y_3,v_0y_0y_2)=1$
and
$$
v_3y_3^2 \eqm{v_0y_1^3y_2^2}{k_4v_2y_0^2}.
$$
Now it  straightforward to deduce from \eqref{ut}, \eqref{capri-5}, 
\eqref{temple}
and the coprimality relation
$\hcf(y_3,v_0y_0y_2)=1$,  that
\begin{align*}
\hcf(v_0y_1^3y_2^2,k_4v_2y_0^2)
&= \hcf(v_0y_1^3y_2^2,k_4)\\
&= \hcf(v_0y_1^3y_2^2, v_1v_2,v_3y_3^2)\\
&= \hcf( \hcf(v_0y_2^2,v_1),v_3y_3^2)) =1,
\end{align*}
for any $k_4$ dividing $v_1v_2$ and $y_4$.
Similarly one sees that
$\hcf(v_3,k_4v_2y_0^2)=1,$
for any such $k_4$.
We shall therefore only be interested in summing over divisors
$k_4\mid v_1v_2$ for which $\hcf(k_4,v_0v_3y_1y_2)=1.$
In fact it suffices to sum over all divisors $k_4\mid v_1v_2$ for
which $\hcf(k_4,v_0v_3y_2)=1,$
since any divisor of $v_1v_2$ is
coprime to $y_1$ by \eqref{capri-5}.
Under this understanding it is now clear
that there exists a unique positive integer
$\vr$, with $\vr\leq k_4v_2y_0^2$ and $\hcf(\vr, k_4v_2y_0^2)=1$, such that
$$
v_3\vr^2 \eqm{v_0y_1}{k_4v_2y_0^2}, \qquad
y_3 \eqm{\vr y_1y_2}{k_4v_2y_0^2}.
$$
Our investigation has therefore led to the equality
$$
S=\sum_{\colt{k_4 \mid 
v_1v_2}{\hcf(k_4,v_0v_3y_2)=1}}\mu(k_4)\sum_{\tolt{\vr\leq 
k_4v_2y_0^2}{v_3\vr^2
      \equiv{v_0y_1}\mod{k_4v_2y_0^2}}{\hcf(\vr,k_4v_2y_0^2)=1}
} S_{k_4}(\vr),
$$
where
$$
S_{k_4}(\vr)=
\#\left\{ y_3\in\Z_*:
\begin{array}{ll}
\hcf(y_3,v_0y_2)=1,~\mbox{\eqref{rolly} holds},\\
y_3 \equiv{\vr y_1y_2}\mod {k_4v_2y_0^2}
\end{array}
\right\}.
$$
Here the coprimality relation $\hcf(y_0,y_3)=1$ follows from
the relations \eqref{capri-5}, \eqref{temple} and $\hcf(\vr,k_4v_2y_0^2)=1$.

In view of the fact that $\hcf(k_4,v_0v_3y_2)=1$,
it follows from \eqref{capri-5} and \eqref{temple} that
$\hcf(\vr y_1y_2,k_4v_2y_0^2)=1$ in the definition of $S_{k_4}(\vr)$.
In order to estimate $S_{k_4}(\vr)$ we may therefore employ Lemma \ref{tek}
with $\kappa=0$ and 
the characteristic function
$$
\chi(n)= \left\{
\begin{array}{ll}
1, & \mbox{if $\hcf(n,v_0y_2)=1$},\\
0, & \mbox{otherwise}.
\end{array}
\right.
$$
Now it is easy to see that
$$
\sum_{\colt{d=1}{\hcf(d,k_4v_2y_0^2)=1}}^\infty\frac{(\chi*\mu)(d)}{d}=
\prod_{\colt{p\mid v_0y_2}{p \nmid k_4v_2y_0}}
\Big(1-\frac{1}{p}\Big)=
\prod_{p\mid v_0y_2} \Big(1-\frac{1}{p}\Big),
$$
whence
$$
S_{k_4}(\vr)=
\phi^*(v_0y_2) \frac{Y_3f(y_1/Y_1,v_1/V_1)}{k_4v_2 y_0^2}
+ O\big(2^{\om(v_0y_2)}\big).
$$
Here, as throughout this paper, we use the notation
\begin{equation}\lab{phi*}
\phi^*(n)=\frac{\phi(n)}{n}=\prod_{p\mid n}\Big(1-\frac{1}{p}\Big)
\end{equation}
for any $n\in \N$.  
Note that the number of positive integers
$\vr\leq k_4v_2y_0^2$ such that $\hcf(\vr, k_4v_2y_0^2)=1$ and
$$
v_3\vr^2 \eqm{v_0y_1}{k_4v_2y_0^2},
$$
is at most $\eta(v_0v_3y_1;k_4v_2y_0^2)\leq 2^{\om(k_4v_2y_0)+1}\leq
2^{\om(v_1v_2y_0)+1}$ by \eqref{triv}.
We have therefore established the
following result.

\begin{lem}\lab{Sumy3y4}
Let $(\ma{v},y_0,y_1,y_2) \in \N^5\times\Z_*\times\N$ satisfy \eqref{capri-5}, \eqref{temple},
\eqref{holly} and \eqref{bolly}.  Then for any $B\geq 1$ we have
$$
S=\frac{Y_3f(y_1/Y_1,v_1/V_1)}{v_2 y_0^2}\Sigma(\ma{v},y_0,y_1,y_2)
   + O\big(2^{\om(v_0y_2)}4^{\om(v_1v_2y_0)}\big),
$$
where
\begin{equation}\lab{sigma1}
\Sigma(\ma{v},y_0,y_1,y_2)=
\phi^*(v_0y_2)
\sum_{\colt{ k_4 \mid v_1v_2}{\hcf(k_4,v_0v_3y_2)=1}}
\hspace{-0.4cm}\frac{\mu(k_4)}{k_4}
\sum_{\tolt{\vr\leq k_4v_2y_0^2}{v_3\vr^2
      \equiv{v_0y_1}\mod{k_4v_2y_0^2}}{\hcf(\vr,k_4v_2y_0^2)=1}}
\hspace{-0.5cm}1.
\end{equation}
\end{lem}

We close this section by showing that once summed over all
$(\ma{v},y_0,y_1,y_2) \in \N^5\times\Z_*\times\N$  satisfying
\eqref{holly} and  \eqref{bolly}, the error term
in Lemma  \ref{Sumy3y4} is satisfactory.
For this we shall make use of the familiar estimate
$$
\sum_{n \leq x}a^{\omega(n)}\ll_a x(\log x)^{a-1},
$$
for any $a\in \N$, in addition to estimates that follow from applying partial
summation to it.
In this way we therefore obtain the overall contribution
\begin{align*}
&\ll \sum_{v_0,v_1,v_2,v_3,y_0,y_2}
\sum_{-Y_1<y_1\leq V_1Y_1/v_1}2^{\om(v_0y_2)}4^{\om(v_1v_2y_0)}\\
&\ll \sum_{v_0,v_1,v_2,v_3,y_0,y_2}
2^{\om(v_0y_2)}4^{\om(v_1v_2y_0)}
\frac{V_1Y_1}{v_1}\\
&\ll (\log B)^4\sum_{v_0,v_2,v_3,y_0,y_2}
2^{\om(v_0y_2)}4^{\om(v_2y_0)}V_1Y_1.
\end{align*}
But now we may employ \eqref{deV1} to bound this quantity by
$$
\ll B^{5/6}(\log B)^4\sum_{v_0,v_2,v_3,y_0,y_2}
\frac{2^{\om(v_0y_2)}4^{\om(v_2y_0)}}{v_0^{7/3}
v_2^{13/6}v_3^{3/2}y_0^{4/3}y_2^{5/3}}
\ll B^{5/6}(\log B)^4.
$$
We shall see below that this is satisfactory.

\subsection{Summation over the variables $y_1$ and $y_2$}\lab{y1}

Our next task is to sum $S$ over all non-zero integers $y_1$ which
satisfy \eqref{capri-5} and \eqref{holly}, and all positive integers
$y_2$ which satisfy $\hcf(y_2,v_2v_3y_0)=1$  and \eqref{golly}.
We therefore write
$$
S'=\frac{Y_3}{v_2 y_0^2}
\sum_{\colt{y_2\leq Y_2}{\hcf(y_2,v_2v_3y_0)=1}}
\sum_{\colt{-Y_1<y_1\leq
V_1Y_1/v_1}{\hcf(y_1,v_0v_1v_2v_3y_0)=1}}
f(y_1/Y_1,v_1/V_1)\Sigma(\ma{v},y_0,y_1,y_2),
$$
where $\Sigma(\ma{v},y_0,y_1,y_2)$ is given by \eqref{sigma1}.

Let $t>0$.
We begin by establishing asymptotic formulae
for the two quantities
$$
\mcal{S}(\pm t)= \phi^*(v_0y_2)
\sum_{\colt{ k_4 \mid v_1v_2}{\hcf(k_4,v_0v_3y_2)=1}}
\frac{\mu(k_4)}{k_4}
\sum_{\colt{\vr\leq k_4v_2y_0^2}{\hcf(\vr,k_4v_2y_0^2)=1}}
S_{k_4}'(\vr;\pm t),
$$
where
\begin{align*}
S_{k_4}'(\vr;+t)
&=\#\left\{ 0\leq y_1 \leq t:
\begin{array}{ll}
\hcf(y_1,v_0v_1v_2v_3y_0)=1, \\
v_3\vr^2 \equiv{v_0y_1}\mod{k_4v_2y_0^2}
\end{array}
\right\},\\
S_{k_4}'(\vr;-t)
&=\#\left\{ -t\leq y_1\leq 0:
\begin{array}{ll}
\hcf(y_1,v_0v_1v_2v_3y_0)=1, \\
v_3\vr^2 \equiv{v_0y_1}\mod{k_4v_2y_0^2}
\end{array}
\right\}.
\end{align*}
Now it is clear that we have
\begin{equation}\lab{hills}
\hcf(v_3\vr^2,k_4v_2y_0^2)=1
\end{equation}
in the definition of $S_{k_4}'(\vr;\pm t)$, since $\hcf(k_4,v_3)=1$.
In particular it follows that we
may replace the coprimality relation
appearing in $S_{k_4}'(\vr;\pm t)$ by $\hcf(y_1,v_0v_1v_3)=1$.
We shall treat this coprimality condition with a M\"obius inversion.
Thus we find that $\mcal{S}(\pm t)$ is equal to
$$
\phi^*(v_0y_2)\hspace{-0.5cm}
\sum_{\colt{ k_4 \mid v_1v_2}{\hcf(k_4,v_0v_3y_2)=1}}\frac{\mu(k_4)}{k_4}
\sum_{\colt{k_1\mid v_0v_1v_3}{\hcf(k_1,k_4v_2y_0)=1}}
\hspace{-0.2cm}
\mu(k_1)
\hspace{-0.5cm}
\sum_{\colt{\vr\leq k_4v_2y_0^2}{\hcf(\vr,k_4v_2y_0^2)=1}}
\hspace{-0.2cm} S_{k_1,k_4}'(\vr;\pm t),
$$
where
\begin{align*}
S_{k_1,k_4}'(\vr;+t)
&=\#\left\{ 0\leq y_1 \leq t/k_1: v_3\vr^2 \equiv{k_1v_0y_1}\mod{k_4v_2y_0^2}
\right\},\\
S_{k_1,k_4}'(\vr;-t)
&=\#\left\{ -t/k_1\leq y_1\leq 0: v_3\vr^2 \equiv{k_1v_0y_1}\mod{k_4v_2y_0^2}
\right\}.
\end{align*}
Here we have used \eqref{hills} to deduce that we must only sum over values
of $k_1\mid v_0v_1v_3$ for which $\hcf(k_1,k_4v_2y_0)=1$.

Let $(\ma{v},y_0) \in \N^5$  satisfy the constraints
\begin{equation}\lab{zolly}
|\mu(v_0v_2v_3)|=\hcf(v_0v_3,y_0 )=1, \quad
v_0^4v_1^6v_2^5v_3^3y_0^4\leq B,
\end{equation}
that follow from \eqref{temple} and \eqref{bolly}. We let
$b_{\pm} \leq k_4v_2y_0^2$ be the unique
positive integer such that
$$
b_{\pm}k_1v_0\equiv{\pm v_3}\mod{k_4v_2y_0^2}.
$$
In particular it follows from \eqref{hills} that 
$\hcf(b_{\pm},k_4v_2y_0^2)=1$, and we
may therefore employ Lemma \ref{cong} to deduce that
$$
S_{k_1,k_4}'(\vr;\pm t) =
\frac{t}{k_1k_4v_2y_0^2}+r(\pm t; b_\pm \vr^2),
$$
where
\begin{equation}\lab{viva}
r(\pm t;b_\pm \vr^2)=\psi\Big( \frac{-b_{\pm}\vr^2}{k_4v_2y_0^2}\Big)
-\psi\Big(\frac{t/k_1-b_{\pm}\vr^2}{k_4v_2y_0^2}\Big).
\end{equation}

Recall the definition \eqref{phi*} of $\phi^*$ and observe that
$
\phi^*(ab)\phi^*(\hcf(a,b))=\phi^*(a)\phi^*(b),
$
for any $a,b\in \N$.
We define
\begin{equation}\lab{vt'}
\vt(\ma{v},y_0,y_2)=\left\{
\begin{array}{ll}
\displaystyle{\frac{ \phi^*(v_0v_1v_2y_2)\phi^*(v_0v_1v_2v_3y_0)
}{\phi^*(\hcf( v_1 ,v_3 ))}}, & \mbox{if \eqref{temple} holds},\\
0, & \mbox{otherwise}.
\end{array}
\right.
\end{equation}
Then a straightforward calculation reveals that
\begin{equation}\lab{ganges}
\mcal{S}(\pm t)= \vt(\ma{v},y_0,y_2) t+ \mcal{R}(\pm t)
\end{equation}
for any non-zero $t>0$,
where
$$
\mcal{R}(\pm t)=\phi^*(v_0y_2)\hspace{-0.5cm}
\sum_{\colt{ k_4 \mid v_1v_2}{\hcf(k_4,v_0v_3y_2)=1}}
\hspace{-0.3cm}
\frac{\mu(k_4)}{k_4}
\hspace{-0.3cm}
\sum_{\colt{k_1\mid v_0v_1v_3}{\hcf(k_1,k_4v_2y_0)=1}}
\hspace{-0.3cm}
\mu(k_1)
\hspace{-0.3cm}
\sum_{\colt{\vr\leq k_4v_2y_0^2}{\hcf(\vr,k_4v_2y_0^2)=1}}
\hspace{-0.3cm} 
r(\pm t;b_\pm \vr^2).
$$
Here $r(\pm t;b_\pm \vr^2)$ is given by \eqref{viva} and the positive integers
$b_-,b_+$ are uniquely determined by fixed choices of
$k_1,k_4,v_0,v_2,v_3, y_0$, as outlined above.

We may now apply partial summation to estimate $S'$.
Now it is clear that $S'=S_-'+S_+',$
where $S_-'$ denotes the contribution from  $y_1$ contained in the
interval $(-Y_1,0]$ and $S_+'$ denotes the contribution from  $y_1$
contained in the
interval $(0,V_1Y_1/v_1]$.  We begin by estimating $S_-'$, for which
we first deduce from  \eqref{range1} and \eqref{range2} that
$$
\frac{v_1}{V_1}=\Big(\frac{y_2}{Y_2}\Big)^{1/3},\quad
Y_1=\Big(\frac{Bv_2y_0^2}{v_0y_2^2}\Big)^{1/3}= \frac{Y_1'}{y_2^{2/3}},
$$
say.  We  may now apply \eqref{ganges}, in conjunction with partial
summation, in order to deduce that
$$
S_-'= \sum_{\colt{y_2\leq Y_2}{\hcf(y_2,v_2v_3y_0)=1}}
\Big(\frac{\vt(\ma{v},y_0,y_2)Y_1Y_3}{v_2y_0^2} 
\int_{-1}^{0}f(u,v_1/V_1)\d u\Big) +
R_-',
$$
where
\begin{align*}
R_-'&=  \frac{Y_3}{v_2 y_0^2}\sum_{\colt{y_2\leq Y_2}{\hcf(y_2,v_2v_3y_0)=1}}
\int_{0}^{1}f'(-u,(y_2/Y_2)^{1/3})\mcal{R}(-uY_1'/y_2^{2/3})\d u\\
&=
\frac{Y_3}{v_2y_0^2}\sum_{\colt{ k_4 \mid v_1v_2}{\hcf(k_4,v_0v_3 )=1}}
\frac{\mu(k_4)}{k_4}
\sum_{\colt{k_1\mid v_0v_1v_3
}{\hcf(k_1,k_4v_2y_0)=1}} \mu(k_1)
\sum_{\colt{\vr\leq k_4v_2y_0^2}{\hcf(\vr,k_4v_2y_0^2)=1}}
\\& \quad
\sum_{\colt{y_2\leq
Y_2}{\hcf
(y_2,k_4v_2v_3y_0)=1}}
\hspace{-0.3cm}
\phi^*(v_0y_2)
\int_{0}^{1}f'(-u,(y_2/Y_2)^{1/3})
r(-uY_1'/y_2^{2/3};b_- \vr^2)\d u.
\end{align*}
Define the arithmetic function
$$
\phi^\dagger(n)=\prod_{p\mid n}\Big(1+\frac{1}{p}\Big)^{-1}.
$$
We estimate $R_-'$ via an application of
Lemma \ref{tek} with 
$a=0,q=1$ and  $\kappa=\ve$.  This gives
\begin{align*}
\sum_{\colt{y_2\leq t}{\hcf (y_2,k_4v_2v_3y_0)=1}}
\hspace{-0.5cm}
\phi^*(v_0y_2)
&=\frac{6}{\pi^2}\phi^\dagger(k_4v_0v_2v_3y_0)t
+O(2^{\omega(v_1v_2v_3y_0)}t^\ve).
\end{align*}
Indeed, the corresponding Dirichlet series is equal to
$$
\phi^*(v_0)\zeta(s)\prod_{ p \nmid
k_4v_0v_2v_3y_0}\Big(1-\frac{1}{p^{s+1}}\Big)
\prod_{p\mid
     k_4v_2v_3y_0} \Big(1-\frac{1}{p^{s }}\Big).
$$
An application of partial
summation therefore yields the estimate
\begin{align}\begin{split}
R_-' 
&= \frac{\varphi_-(\ma{v},y_0)Y_2Y_3}{v_2y_0^2}+
O \Big(
{2^{\omega(v_1v_2)+\omega(v_0v_1v_3)+\omega(v_1v_2v_3y_0)}}Y_2^\ve Y_3\Big),\\
&= \frac{\varphi_-(\ma{v},y_0)Y_2Y_3}{v_2y_0^2}+
O_\ve \big( B^\ve Y_3\big),\end{split}\lab{estR-'}
\end{align}
where
\begin{align*}
\varphi_-(\ma{v},y_0)&=\frac{18}{\pi^2}
\hspace{-0.3cm}
\sum_{\colt{ k_4 \mid v_1v_2}{\hcf(k_4,v_0v_3 )=1}}
\hspace{-0.3cm}
\frac{\mu(k_4)\phi^\dagger(k_4v_0v_2v_3y_0)}{k_4}
\hspace{-0.3cm}
\sum_{\colt{k_1\mid v_0v_1v_3
}{\hcf(k_1,k_4v_2y_0^2)=1}}
\hspace{-0.3cm}
\mu(k_1)
\\&\quad
\int_0^1\int_0^1
\Big(t^2f'(-u,t)
\hspace{-0.3cm}
\sum_{\colt{\vr\leq k_4v_2y_0^2}{\hcf(\vr,k_4v_2y_0^2)=1}}
\hspace{-0.3cm}
r(-v_0v_1^2v_2^2v_3y_0^2u/t^2;b_- \vr^2)\Big)\d u \d t.
\end{align*}
Here we have used the trivial inequality $2^{\omega(n)}=O_\ve(n^\ve)$
for any $n \in \N$.
An application of Lemma \ref{av-order'} clearly reveals that
$$
\varphi_-(\ma{v},y_0) \ll_\ve 
(v_2y_0^2)^{1/2+\ve}2^{\omega(v_1v_2)+\omega(v_0v_1v_3)},
$$
for any $\ve>0$. 
Our estimate \eqref{estR-'} for  $R_-'$ isn't terribly good when $Y_2$ is
small. Fortunately, by inverting the order of summation over $\vr$ and $y_2$
we may use Lemma \ref{av-order'} to deduce the alternative  
estimate 
\begin{align}\begin{split}
R_-' &= \frac{\varphi_-(\ma{v},y_0)Y_2Y_3}{v_2y_0^2} +O_\ve\Big(
{2^{\omega(v_1v_2)+\omega(v_0v_1v_3) }}
\frac{Y_2Y_3}{(v_2y_0^2)^{1/2-\varepsilon}}\Big)\\ 
&= \frac{\varphi_-(\ma{v},y_0)Y_2Y_3}{v_2y_0^2} +O_\ve\Big(
B^{\varepsilon}\frac{Y_2Y_3}{(v_2y_0^2)^{1/2}}\Big).
\end{split}\label{R-'compl}
\end{align}
Note here that the main term is dominated by the error term.
On combining \eqref{estR-'} and \eqref{R-'compl}, however, we obtain
the estimate
\begin{align*}
R_-'  =& \frac{\varphi_-(\ma{v},y_0)Y_2Y_3}{v_2y_0^2} +
O_\ve\Big( B^{\varepsilon}
Y_3\min\Big\{1,\frac{ Y_2}{(v_2y_0^2)^{1/2 }}\Big\}\Big). 
\end{align*}

Arguing in a similar fashion it is straightforward to deduce that
$$
S_+'= \sum_{\colt{y_2\leq Y_2}{\hcf(y_2,v_2v_3y_0)=1}}
\Big(\frac{\vt(\ma{v},y_0,y_2)Y_1Y_3}{v_2y_0^2} 
\int_{0}^{V_1/v_1}f(u,v_1/V_1)\d u\Big) +
R_+',
$$
where
\begin{align*}
R_+'  =& \frac{\varphi_+(\ma{v},y_0)Y_2Y_3}{v_2y_0^2}+O_\ve\Big( 
B^{\varepsilon}
Y_3\min\Big\{1,\frac{ Y_2}{(v_2y_0^2)^{1/2 }}\Big\}\Big).
\end{align*}
Here one finds that 
\begin{align*}
\varphi_+(\ma{v},y_0)&=\frac{18}{\pi^2}
\hspace{-0.3cm}
\sum_{\colt{ k_4 \mid v_1v_2}{\hcf(k_4,v_0v_3 )=1}}
\hspace{-0.3cm}
\frac{\mu(k_4)\phi^\dagger(k_4v_0v_2v_3y_0)}{k_4}
\hspace{-0.3cm}
\sum_{\colt{k_1\mid v_0v_1v_3
}{\hcf(k_1,k_4v_2y_0^2)=1}}
\hspace{-0.3cm}
\mu(k_1)
\\&\quad
\int_0^1\int_0^{1/t}
\Big(t^2f'(u,t)
\hspace{-0.3cm}
\sum_{\colt{\vr\leq k_4v_2y_0^2}{\hcf(\vr,k_4v_2y_0^2)=1}}
\hspace{-0.3cm}
r(v_0v_1^2v_2^2v_3y_0^2u/t^2;b_+ \vr^2)\Big)\d u \d t,
\end{align*}
with $\varphi_+(\ma{v},y_0) \ll_\ve
(v_2y_0^2)^{1/2+\ve}2^{\omega(v_1v_2)+\omega(v_0v_1v_3)}$.

We may now complete our estimate for $S'$.
Recall the definition \eqref{f} of the function $f(u,v)$, and define
\begin{equation}\lab{g}
g(v)=\int_{-1}^{1/v}f(u,v)\d u.
\end{equation}
Then $g$ is a bounded differentiable function, whose derivative is
also bounded on the interval $[0,\infty)$.
Moreover let
\begin{equation}\lab{vpi}
\varphi(\ma{v},y_0)= \varphi_-(\ma{v},y_0)+\varphi_+(\ma{v},y_0).
\end{equation}
Then on combining our various estimates we have therefore established
the following result.

\begin{lem}\lab{Sumy1}
Let $(\ma{v},y_0) \in \N^5$  satisfy \eqref{zolly}.
Then for any $B\geq 1$ we have
\begin{align*}
S'=
  \sum_{\colt{y_2\leq Y_2}{\hcf(y_2,v_2v_3y_0)=1}}
& \Big(\frac{\vt(\ma{v},y_0,y_2)Y_1Y_3g(v_1/V_1)}{v_2 y_0^2}\Big) 
+\frac{\varphi(\ma{v},y_0)Y_2Y_3}{v_2y_0^2}\\
&+O_\ve\Big( B^{\varepsilon}
Y_3\min\Big\{1,\frac{ Y_2}{(v_2y_0^2)^{1/2 }}\Big\}\Big), 
\end{align*}
where $\vt(\ma{v},y_0,y_2)$ is given by \eqref{vt'}, $g$ is given by
\eqref{g} and $\varphi(\ma{v},y_0)$ is given by \eqref{vpi} and
satisfies
\begin{equation}\lab{uea}
\varphi(\ma{v},y_0) \ll_\ve 
(v_2y_0^2)^{1/2+\ve}2^{\omega(v_1v_2)+\omega(v_0v_1v_3)},
\end{equation}
for any $\ve>0$.
\end{lem}

We end this section by
showing that once summed over all
$(\ma{v},y_0) \in \N^5$  satisfying
\eqref{zolly}, the error term
in Lemma  \ref{Sumy1} is satisfactory.
On recalling the definition \eqref{range2} of $Y_2$ and $Y_3$, and then first summing
over $y_0$, we easily obtain the satisfactory overall contribution
\begin{align*}
&\ll_\ve
B^{1/2+\varepsilon}\sum_{ v_0,v_1,v_2,v_3, y_0}
\frac{(v_2y_0^2)^{1/2}}{v_3^{1/2}}\min\Big\{1,
\frac{B^{1/2}}{v_0^2v_1^3v_2^3v_3^{3/2}y_0^3}\Big\}  
\\ &\ll_\ve B^{5/6+\ve} 
\sum_{ v_0,v_1,v_2,v_3}
\frac{1}{v_0^{4/3}v_1^{2}v_2^{3/2}v_3^{3/2}}\ll_\ve
B^{5/6+ \ve}.
\end{align*}

\subsection{Summation over the remaining variables}\lab{final}

In this section we complete our preliminary estimate for $N_{U,H}(B)$.
It is clear from Lemma \ref{Sumy1} that we have two distinct terms to
deal with.
We begin by deducing from \eqref{range2} that
$$
\frac{Y_1Y_3}{v_2y_0^2}=\frac{B^{5/6} n^{1/6}}{ v_0v_1v_2 v_3
y_0 y_2},
$$
in the statement of Lemma \ref{Sumy1}, with $n=v_0^4v_1^6
v_2^5v_3^3y_0^4y_2^2$.
Define the arithmetic function
\begin{equation}\lab{46}
\Delta(n)=B^{-5/6}\sum_{\colt{\ma{v},y_0,y_2}
{v_0^4v_1^6 v_2^5v_3^3y_0^4y_2^2=n}}\frac{\vt(\ma{v},y_0,y_2)Y_1Y_3}{v_2y_0^2},
\end{equation}
where $\vt(\ma{v},y_0,y_2)$ is given by \eqref{vt'}. 
Recall the definition 
\eqref{g} of the function $g$ and that of the counting function $N(Q_1,Q_2;B)$ that
appears in the statement of Lemma \ref{red-1}.
Let $\ve>0$.
We proceed by 
establishing the existence of a constant
$\be \in \R$ for which 
\begin{equation}\lab{cof}
N (Q_1,Q_2;B )= B^{5/6}\sum_{n\leq
B}\Delta(n)g \Big(\Big(\frac{n}{B}\Big)^{1/6}\Big)+ \be B + O_\ve(B^{5/6+\ve}).
\end{equation}
This follows rather easily from Lemma \ref{Sumy1}.  
Define the sum
$$
T(B)=\sum_{\colt{\ma{v},y_0}{\mbox{\scriptsize{\eqref{zolly} holds}}}}
\frac{\varphi(\ma{v},y_0)Y_2Y_3}{v_2y_0^2},
$$
for any $B \geq 1$.  Then in view of the error terms that we have
estimated along the way in \S \ref{y3y4} and \S \ref{y1}, 
it is clearly enough to 
establish the existence of a constant
$\be \in \R$ for which
$$
T(B)= \be B + O_\ve(B^{5/6+\ve}).
$$
On recalling \eqref{range2}, we see that
$$
\frac{Y_2Y_3}{v_2y_0^2}=\frac{B}{v_0^2v_1^3v_2^3v_3^2y_0^3}.
$$
Hence on taking $\varepsilon<1/3$, it follows from \eqref{uea} that 
\begin{align*}
T(B)-\be B
&\ll_\ve B\sum_{\colt{\ma{v},y_0}{v_0^4v_1^6 v_2^5v_3^3y_0^4> B}}
\frac{(v_0v_1v_2v_3y_0)^\ve}{v_0^2v_1^3v_2^{5/2}v_3^2y_0^2}\\
&\ll_\ve B^{5/6} \sum_{\ma{v},y_0}
\frac{(v_0v_1v_2v_3y_0)^\ve}{v_0^{4/3}v_1^{2}v_2^{5/3}v_3^{3/2}y_0^{4/3}}
\ll_\ve B^{5/6},
\end{align*}
with
\begin{equation}\lab{beta'}
\be=\sum_{\colt{\ma{v},y_0}{\hcf(v_0v_3,y_0)=1}}
\frac{|\mu(v_0v_2v_3)|\varphi(\ma{v},y_0)}{v_0^2v_1^3v_2^3v_3^2y_0^3}.
\end{equation}
This therefore completes the proof of \eqref{cof}.
On inserting this estimate into Lemma \ref{red-1} we obtain the following result.

\begin{lem}\lab{Sum-all'}
Let $\ve>0$.  Then for any $B \geq 1$ we have
$$
N_{U,H}(B)= 2B^{5/6}\sum_{n\leq
B}\Delta(n)g \Big(\Big(\frac{n}{B}\Big)^{1/6}\Big)+ \Big( \frac{12}{\pi^2}+2\be\Big)B + O_\ve(B^{5/6+\ve}),
$$
where $g$ is given by \eqref{g}, $\D$ is given by 
\eqref{46} and $\beta$ is given by \eqref{beta'}.
\end{lem}

\section{The height zeta function}\lab{height}

For $\Re e(s)>1$ we recall the definition of the 
height zeta function \eqref{h-zeta}, and the identity \eqref{integral}.
Thus it follows from Lemma \ref{Sum-all'} that $Z_{U,H}(s)=Z_1(s)+Z_2(s)$, where 
\begin{align*}
Z_1(s)&=
2s\int_1^\infty t^{-s-1/6}\sum_{n\leq t}\Delta(n)g
\Big(\Big(\frac{n}{t}\Big)^{1/6}\Big)\d t,\\
Z_2(s)&=\frac{{12/\pi^2}+2\be}{s-1}+G_2(s),
\end{align*} 
and 
$$
G_2(s)=s\int_{1}^\infty t^{-s-1}R(t)\d t 
$$ 
for some function $R(t)$ such that $R(t)\ll_\varepsilon
t^{5/6+\varepsilon}$ for any $\varepsilon>0$.  But then it easily
follows that 
$G_2(s)$ is holomorphic on the half-plane $\Re e(s)\geq 5/6+\varepsilon,$ and satisfies 
$G_2(s)\ll 1+|\Im m (s)|$
on this domain.  Finally an application of the  Phragm\'en-Lindel\"of
Theorem yields the finer upper bound
$$
G_2(s)\ll_\ve (1+|\Im m (s)|)^{6(1-\Re e (s))+\varepsilon}.
$$
on this domain.

To establish Theorem \ref{main'} it therefore remains to
analyse the function $Z_1(s)$.
Recall the definition \eqref{46} of $\D$ and define the corresponding Dirichlet series
$$
F(s)=\sum_{n=1}^\infty \frac{\Delta(n)}{n^s}.
$$
Then it is easily seen
that 
$$
Z_1(s)=2sF(s-5/6)
\int_1^\infty t^{-s-1/6}g(1/t^{1/6})\d t
=F(s-5/6)G_{1,1}(s),
$$
where 
\begin{equation}
\label{G1}
G_{1,1}(s) =12s\int_0^1 v^{6s-6}g(v)\d v.
\end{equation}
Recall the definition \eqref{g} of $g$.
Then a simple calculation reveals that $G_{1,1}(1) =12\tau_\infty$, in the notation of \eqref{t-inf}.
Moreover, an application of partial integration yields
$$
G_{1,1}(s) =\frac{12s}{6s-5}\Big(g(1)-\int_0^1
v^{6s-5}g'(v)\d v\Big),
$$
whence it is clear that $G_{1,1}(s)$ 
is holomorphic and bounded
on the half-plane $\Re e(s)\geq 5/6+\ve$ for any $\ve>0$.

We proceed by analysing the Dirichlet series $F(s-5/6)$ in more detail.
Define the function 
\begin{equation}\lab{G2}
G_{1,2}(s)= \frac{F(s-5/6)}{E_1(s)E_2(s)},
\end{equation}
for $\Re e(s)>5/6$ and let $\ve>0$.  
Here $E_1(s)$ and $E_2(s)$ are given by \eqref{e1} and \eqref{e2}, respectively.
In order to complete the proof of
Theorem \ref{main'}, with
\begin{equation}\lab{G!}
G_1(s)=G_{1,1}(s)G_{1,2}(s),
\end{equation}
it remains to establish that 
$G_{1,2}(1) \neq 0$ and that 
$G_{1,2}(s)$ is holomorphic and bounded for  $\Re
e(s)\geq 5/6+\varepsilon$.   This is achieved for us in the following result.

\begin{lem}\label{balloon}
Let $\varepsilon>0$. Then $G_{1,2}(s+1)$ 
is holomorphic and bounded on the half-plane 
$\mcal{H}=\{ s\in \C  : \Re e(s)\geq -1/6+\varepsilon\}$.
\end{lem}

\begin{proof}
On writing 
$$
G_{1,2}(s+1)=\prod_p G_{p}(s+1),
$$
it will clearly suffice to show that
$G_{p}(s+1)= 1+O_\ve(1/p^{1+\varepsilon})$ uniformly on $\mcal{H}$.
We begin the proof of Lemma \ref{balloon} by observing that 
$$
F(s+1/6)=
\hspace{-0.3cm}
\sum_{\tolt{(\ma{v},y_0,y_2)\in 
\N^6}
{\hcf(v_2v_3y_0,y_2)=1}{\hcf(v_0v_3,y_0)=1}} 
\hspace{-0.3cm}
\frac{|\mu(v_0v_2v_3)|
\phi^*(v_0v_1v_2y_2)\phi^*(v_0v_1v_2v_3y_0)}{\phi^*(\hcf( 
v_1 ,v_3)
)v_0^{4s+1}v_1^{6s+1} 
v_2^{5s+1}v_3^{3s+1}y_0^{4s+1}y_2^{2s+1}}.
$$
Thus on writing $F(s+1/6) = \prod_p F_p(s+1/6)$ as a product of
local factors, a straightforward calculation reveals that
$F_{p}(s+1/6)$ is equal
to
\begin{equation}\lab{saw}
\begin{split}
1 &+\frac{1-1/p}{p^{2s+1}-1}+\frac{1-1/p}{p^{4s+1}-1}
+\frac{(1-1/p)^2}{ 
p^{6s+1}-1}\Big(\frac{p^{2s+1}}{p^{2s+1}-1}+\frac{1}{p^{4s+1}-1}\Big)\\
&+ \frac{p^{4s+1}(1-1/p)^2}{ (p^{2s+1}-1)(p^{6s+1}-1)}
+ \frac{p^{5s+1}(1-1/p)^2}{ (p^{4s+1}-1)(p^{6s+1}-1)} 
+ \frac{p^{3s}(1-1/p)}{ p^{6s+1}-1},
\end{split}
\end{equation}
for any prime $p$. 
On collecting together factors of $(p^{2s+1}-1)^{-1}$ and
$(p^{4s+1}-1)^{-1}$ we see that
\begin{align*}
F_p(s+1/6)\Big(1-\frac{1}{p^{6s+1}}\Big)
&=1- \frac{1}{p^{6s+1}} +\frac{1}{p^{3s+1}}\\
&\quad+\frac{1}{p^{2s+1}-1}\Big(1 +
\frac{1}{p^{2s}} +\frac{1}{p^{4s}} 
-\frac{1}{p^{6s+1}}\Big)\\
&\quad +\frac{1}{p^{4s+1}-1}\Big(1 +\frac{1}{p^{s}} \Big)
+O_\ve\Big(\frac{1}{p^{1+\varepsilon}}\Big),
\end{align*}
on $\mcal{H}$.  
We now record the obvious estimates 
\begin{align*}
\frac{1}{p^{2s+1}-1}&=\frac{1}{p^{2s+1}
}+\frac{1}{p^{4s+2}}+ O\Big(\frac{1}{p^{2+6\varepsilon}}\Big),
\\\frac{1}{p^{4s+1}-1}&=\frac{1}{p^{4s+1}
}+\frac{1}{p^{8s+2}}+\frac{1}{p^{12s+3}
}+O_\ve\Big(\frac{1}{p^{4/3+16\varepsilon}}\Big),
\end{align*}
and 
$$
1 +\frac{1}{p^{2s}} +\frac{1}{p^{4s}}
-\frac{1}{p^{6s+1}}\ll_\ve p^{2/3-4\varepsilon}, \quad
1 +\frac{1}{p^{s}} \ll_\ve p^{1/6-\varepsilon},
$$
that all hold on $\mcal{H}$.  Combining these estimates we 
therefore deduce that
\begin{align*}
F_p(s+1/6)\Big(1-\frac{1}{p^{6s+1}}\Big)&=1+\frac{1}{p^{2s+1}}+\frac{1}{p^{3s+1}}+\frac{2}{p^{4s+1}}
+\frac{1}{p^{5s+1}}
\\&\quad+\frac{1}{p^{8s+2}}+\frac{1}{p^{9s+2}}+
\frac{1}{p^{13s+3}}+O_\ve\Big(
\frac{1}{p^{1+\varepsilon}}\Big).
\end{align*}
Write $E_{1,p}(s+1)$ for the Euler factor of \eqref{e1}  and write
$E_{2,p}(s+1)$ for the Euler factor of \eqref{e2}.
Then it is now a routine matter to deduce that 
\begin{align*}
\frac{F_p(s+1/6)}{E_{1,p}(s+1)}&=1-\frac{3}{p^{7s+2}}-\frac{ 3}{p^{8s+2}}-\frac{1}{p^{9s+2}}
-\frac{ 1}{p^{10s+2}} +\frac{3}{p^{13s+3}}\\
&\quad+\frac{1}{p^{14s+3}}+O\Big(
\frac{1}{p^{1+\varepsilon}}\Big)\\
&=E_{2,p}(s+1)\Big(1+O\Big(
\frac{1}{p^{1+\varepsilon}}\Big)\Big),
\end{align*}
on $\mcal{H}$, which therefore completes the proof of Lemma \ref{balloon}.
\end{proof}

It remains to combine the expression \eqref{saw} for
$F_p(s+1/6)$, with \eqref{e1} and \eqref{G2} in order to deduce that
$$
E_2(1)G_{1,2}(1)=
\prod_{p}\Big(1-\frac{1}{p}\Big)^6\Big(1+\frac{6
}{p}+\frac{1}{p^2}
\Big) \neq 0.
$$
This therefore completes the proof of Theorem \ref{main'}.

\section{Deduction of Theorem \ref{main}}\lab{deduc}

In this section we shall deduce Theorem \ref{main} from Theorem
\ref{main'} and Lemma \ref{Sum-all'}.  
Let $\ve>0$ and let $T \in [1,B]$.  
Then an application of Perron's formula yields
\begin{align}\begin{split}\lab{pint}
N_{U,H}(B)-
\Big(\frac{12}{\pi^2}+2\be\Big) B =& \frac{1}{2\pi i}
\int_{1+\ve-iT}^{1+\ve+iT} E_1(s)E_2(s)G_1(s)\frac{B^s}{s}\d s \\&+
O_\ve\Big(
\frac{B^{11/6+\ve}}{T}\Big).
\end{split}\end{align}
We apply Cauchy's residue theorem to the rectangular contour $\mcal{C}$ 
joining the points ${\kappa-iT}$, ${\kappa+iT}$,
${1+\ve+iT}$ and ${1+\ve-iT}$, for any $\kappa \in [11/12,1)$.
We must calculate the residue of $E_1(s)E_2(s)G_1(s)B^s/s$ at $s=1$. 
For $\Re e (s)>9/10$ Theorem \ref{main'} implies that the product $E_2(s)G_1(s)$ is holomorphic and bounded.  
In view of \eqref{e1}, we see that
$$ 
E_1(s)=\frac{1}{2880(s-1)^{6}}+O\Big(\frac{1}{(s-1)^{5}}\Big),
$$
as $s\rightarrow 1$.  Hence it follows that 
$$
\rom{Res}_{s=1}\Big\{E_1(s)E_2(s)G_1(s)\frac{B^s}{s}\Big\} =
\frac{E_2(1)G_1(1)}{5!2880}BQ_1(\log B),
$$
for some monic polynomial $Q_1$ of degree $5$. 
Recall from \eqref{G!} that $G_1=G_{1,1}G_{1,2}$.  Then we have already
seen in the previous section that $G_1(1)=12 \tau_\infty \tau$, in the
notation of \eqref{t-inf} and \eqref{deftauf}.
%On recalling the definition \eqref{g} of $g$, it clearly follows from \eqref{t-inf} and 
%\eqref{G1'} that 
%$$
%G_{1,1}(1)=12\tau_\infty.
%$$
%Next we combine the expression \eqref{saw} for
%$F_p(s+1/6)$, with the equation \eqref{G2} and the definition
%\eqref{e1} of $E_1(s)$, in order to deduce that
%$$
%E_2(1)G_{1,2}(1)=
%\prod_{p}\Big(1-\frac{1}{p}\Big)^6\Big(1+\frac{6
%}{p}+\frac{1}{p^2}
%\Big) = \tau,
%$$
%in the notation of \eqref{deftauf}.
Putting all of this together we have therefore shown that
$$
\frac{1}{2\pi i} \int_{\mcal{C}}
E_1(s)E_2(s)G_1(s)\frac{B^s}{s}\d s =
\frac{\tau\tau_\infty}{28800}BQ_2(\log B),  
$$
for some monic polynomial $Q_2$ of degree $5$.
Define the difference
$$
E(B)=N_{U,H}(B)-\frac{\tau\tau_\infty}{28800}BQ_2(\log B)-
\Big(\frac{12}{\pi^2}+2\be\Big) B, 
$$
Then, in view of \eqref{pint} and the fact that the product $E_2(s)G_1(s)$
is holomorphic and bounded for $\Re e (s)>9/10$, we deduce that
\begin{equation}
\begin{split}
E(B) 
&\ll_\ve
\frac{B^{11/6+\ve}}{T}+
\Big(\int_{\kappa-iT}^{\kappa+iT}+\int_{\kappa-iT}^{1+\ve-iT}+ \int_{1+\ve+iT}^{\kappa+iT}\Big)
\Big|E_1(s)\frac{B^s}{s}\Big|\d s,
\lab{t3}
\end{split}
\end{equation}
for any $\kappa \in [11/12,1)$ and any $T\in [1,B]$.
We begin by estimating the contribution from the horizontal
contours. Recall the well-known convexity bound
$$
\zeta(\sigma+i t)\ll_\ve |t|^{(1-\sigma)/3+\varepsilon},
$$
that is valid for any $\sigma\in[1/2,1]$ and $|t|\geq 1$.  
Then it follows that
\begin{equation}\lab{majE1}
E_1(\sigma+it)\ll_\varepsilon
|t|^{8(1-\sigma)+\varepsilon}
\end{equation} 
for any $\sigma\in[11/12,1]$ and $|t|\geq 1$. 
This estimate allows us to deduce that
\begin{equation}\begin{split}
\int_{\kappa-iT}^{1+\ve-iT}\Big|E_1(s)\frac{B^s}{s}\Big|\d s &\ll_\ve 
\int_{\kappa}^{1+\ve}B^\sigma T^{7-8\sigma+\ve}\d \sigma\\ 
&\ll_\ve \frac{B^{1+\ve}T^\ve}{T} + B^{\kappa}T^{7-8\kappa+\ve}.
\lab{t2}
\end{split}
\end{equation}
One obtains the same estimate for the contribution from the remaining
horizontal contour joining ${\kappa+iT}$ to ${1+\ve+iT}$.
\goodbreak

We now turn to the size of the integral
\begin{equation}\lab{t1}
\int_{\kappa-iT}^{\kappa+iT}\Big|E_1(s)\frac{B^s}{s}\Big|\d s \ll  
B^\kappa\int_{-T}^{T} \frac{|E_1(\kappa+it)|}{1+|t|}  \d t =
B^\kappa I(T),
\end{equation}
say.  For given $0<U \ll T$, we begin by estimating the contribution to $I(T)$ from each
integral
$$
\int_{U}^{2U} \frac{|E_1(\kappa+it)|}{1+|t|}\d t \ll 
\frac{1}{U}\int_{U}^{2U}|E_1(\kappa+it)|\d t = \frac{J(U)}{U},
$$
say.  
Let $\ve>0$ and let $k \in \N$.  Then we define $\sigma_k$ to be the infimum of $\sigma$ such that
$$
\frac{1}{T} \int_{1}^T|\zeta(\sigma +i t)|^{2k} \d t = O_\ve(T^{\ve}).
$$
It therefore follows from the mean-value theorem in \cite[\S
7.8]{titch} that
\begin{equation}\lab{mean}
\int_{U}^{2U}|\zeta(\sigma+i t)|^{2k} \d t \ll_\ve U^{1+\ve},
\end{equation}
for any $\sigma \in (\sigma_k,1]$, and any $U \geq 1$.
We shall apply this estimate in the cases $k=2$ and $k=4$, for which
we combine a result due to Heath-Brown \cite{hb} with well-known
estimates for the fourth moment of $|\zeta(1/2+i t)|$ in order to
deduce that
\begin{equation}\lab{dollar}
\sigma_k \leq 
\left\{
\begin{array}{ll}
1/2, & k=2,\\
5/8, & k=4.
\end{array}
\right.
\end{equation}
Returning to our estimate for $J(U)$, for fixed $0<U \ll T$ and any
$\kappa \in [11/12,1)$, we define
$J(U;c)=\int_{U}^{2U}|\zeta(c\kappa-c+1 +cit)|^{4}\d t$. Then we may 
apply H\"older's inequality to deduce that
$$
J(U) \leq J(U;6)^{1/4}J(U;5)^{1/4}J(U;4)^{1/4}J(U;3)^{1/8}J(U;2)^{1/8}.
$$
On combining \eqref{mean}, \eqref{dollar} and the fact that
$\kappa \in[11/12,1)$, we therefore deduce that
$J(U)\ll_\ve U^{1+\ve},$
on re-defining the choice of $\ve$.
Summing over dyadic
intervals for $0<U \ll T$ we obtain
$$
\int_{0}^{T}\frac{|E_1(\kappa+it)|}{1+|t|}\d t \ll_\ve T^\ve.
$$
We obtain the same estimate for the integral over the interval
$[-T,0]$, and so it follows that $I(T)\ll_\ve T^\ve$.  We may
insert this estimate into \eqref{t1}, and then combine it with
\eqref{t2} in \eqref{t3}, in order to conclude that
$$
E(B)\ll_\ve \frac{B^{11/6+\ve}T^\ve}{T}+
B^{\kappa}T^\ve, 
$$
for any $T \in [1,B]$.
We therefore complete the proof of Theorem \ref{main} by taking 
$T=B$.

%On peut am\'eliorer ce r\'esultat en d\'ecalant jusqu'\`a $\kappa=9/10$
%et en appliquant l'in\'egalit\'e d'H\"older avec les exposants
%$(1/4,1/4,3/10,1/10,1/10)$. Pour majorer le premier moment d'ordre $4$ sur
%la droite d'abscisse $2/5$, on applique d'abord l'\'equation
%fonctionnelle de
%$\zeta$ pour se ramener \`a une abscisse $\geq 1/2$. Un calcul analogue
%\`a celui r\'edig\'e ici fournit le th\'eor\`eme \ref{main} avec
%$\delta\in (0,1/11).$

\end{document}